\documentclass[12pt]{article}
\usepackage{amsmath}
\usepackage{amstext}
\usepackage{amsfonts}
\usepackage{amsthm}
\usepackage{amscd}

\numberwithin{equation}{section}

\swapnumbers
\theoremstyle{plain}
\newtheorem{thm}{Theorem}[section]
\newtheorem{lem}[thm]{Lemma}
\newtheorem{prop}[thm]{Proposition}
\newtheorem{cor}[thm]{Corollary}
\theoremstyle{remark}
\newtheorem{rem}[thm]{Remark}

\def\mod0{{\mathcal SU}_C(2,{\mathcal O})}
\def\modk{{\mathcal SU}_C(2,K)}

\def\notin{\in \hspace{-4mm}/\ }

\def\gg{\Gamma_{00}}
\def\ggd{\Gamma_{00}^{(2)}}
\def\ggg{\Gamma_{000}}
\def\ggf{\Gamma_{11}}
\def\cmc{C_2-C_2}
\def\ctc{C_2 \times C_2}

\def\map#1{\ \smash{\mathop{\longrightarrow}\limits^{#1}}\ }
\newcommand{\thetachar}[2]{\genfrac{[}{]}{0pt}{}{#1}{#2}}

\def\sing{\mathrm{Sing} \: \Theta}
\def\spansing{\langle \mathrm{Sing} \Theta \rangle}
\def\tankum{\mathbb{T}_0}
\def\tanmod{\mathbb{T}}
\def\si{\mathrm{Sing} \:}

\def\codim{\mathrm{codim}\:}
\def\sym{\mathrm{Sym}}
\def\ker{\mathrm{ker}\:}
\def\coker{\mathrm{coker}\:}
\def\pic{\mathrm{Pic}}
\def\mult{\mathrm{mult}}
\def\kum{\mathrm{Kum}}
\def\dim{\mathrm{dim} \:}
\def\det{\mathrm{det} \:}
\def\rk{\mathrm{rk} \:}
\def\deg{\mathrm{deg} \:}
\def\Pf{\mathrm{Pf} \:}

\def\lra{\longrightarrow}

\def\ra{\rightarrow}

\def\lms{\longmapsto}

\def\cc{\mathbb{C}}
\def\pp{\mathbb{P}}
\def\ca{\mathcal{C}}

\begin{document}

\title{Singularities of 2$\theta$-divisors in the Jacobian}
\author{C. Pauly \and E. Previato}
\date{}

\maketitle

\section{Introduction}
Let $C$ be a smooth, connected, projective 
{non-hyperelliptic}  
curve of genus $g \geq 2$
over the complex numbers 
and let $\pic^d(C)$ be the
connected component of its Picard variety parametrizing degree $d$
line bundles, for  $d \in \mathbb{Z}$. The variety $\pic^{g-1}(C)$
carries a naturally defined divisor, the Riemann theta divisor $\Theta$,
whose support consists of line bundles that have
nonzero global sections. The
Riemann singularity theorem describes the singular locus of the
divisor $\Theta$ as
$$ \sing = \{ \xi \in \pic^{g-1}(C) \: \mid \: h^0(C,\xi)>1 \}. $$

The Jacobian variety of $C$, $JC := \pic^0(C)$, 
carries a canonically defined line bundle, which we denote by
$\mathcal{O}(2\Theta_0)$. This line bundle is obtained by translating the
bundle $\mathcal{O}(2 \Theta)$ by any theta characteristic. There
is a canonical duality, called the Wirtinger duality, between the
spaces of global sections (\cite{mumf} p. 335):
\begin{equation} \label{wirtinger}
w: |2\Theta|^{\ast} \cong |2\Theta_0|
\end{equation}

In \cite{geem-geer} the authors 
revisit and extend the work of Frobenius on
the subseries $\mathbb{P} \gg \subset |2\Theta_0|$ consisting of
$2\theta$-divisors having multiplicity at least  4 at the origin and
formulate a number of Schottky-type conjectures, some of which have
been proved (\cite{welt} and  \cite{izadi}). Moreover,
J. Fay observed that $2\theta$-divisors
of the Jacobian satisfy a remarkable equivalence (e.g. \cite{welt}
Prop. 4.8, \cite{gu1} Cor. 1)
\begin{equation} \label{frobenius}
\mult_0(D) \geq 4 \iff C-C \subset D \hspace{1cm} \forall D \in |2\Theta_0|
\end{equation}
where the surface $C-C$ denotes the image of the difference map
$\phi_1 : C \times C \ra JC$, which sends a pair of points $(p,q)$ to the
line bundle $\mathcal{O}(p-q)$. We denote by $\mathbb{P}(\mathbb{T}_0)$
the embedded tangent space
at the origin to the Kummer variety $\kum := JC / \pm \subset |2\Theta|$.
It is well-known that the space $\mathbb{T}_0^{\bot}$ of hyperplanes
containing $\mathbb{T}_0$ corresponds under the Wirtinger duality
\eqref{wirtinger} to $\gg$. It follows that
\begin{equation} \label{codimgg}
\codim \gg = 1 + \binom{g+1}{2}.
\end{equation}

\bigskip

In this paper we are interested in some subseries of $\gg$, which
have appeared in the work of several authors (e.g. \cite{gu1},
\cite{gu2}, \cite{welt}). We consider the subseries of $2\theta$-divisors
which are singular along the surface $C-C$, i.e.
\begin{equation} \label{ggd}
\ggd = \{ D \in \gg \: | \: \mult_{p-q}(D) \geq 2 \: \: \forall p,q \in C \} .
\end{equation}

By \eqref{wirtinger}
the morphism given by the linear series $|2\Theta|$ maps $\pic^{g-1}(C)$
to $|2\Theta_0|$. 
Let $\spansing \subset H^0(JC,\mathcal{O}(2 \Theta_0))$ 
denote the linear span of the image of $\sing$ under this morphism. 
We prove:

\vspace{1cm}

\begin{thm}
For any non-hyperelliptic curve
\begin{enumerate}
\item[(1)] there exists a canonical isomorphism 
$$\gg/\ggd \map{\sim} \Lambda^3 H^0(K) $$
\item[(2)] we have an equality 
$$\spansing = \ggd .$$ 
\end{enumerate}
\end{thm}

The method used in the proof of this theorem (section 4) has been
developed in a recent paper by van Geemen and Izadi \cite{geem-iz} and
the key point are the incidence relations (section 2.4) between two
families of stable rank $2$ vector bundles with fixed trivial (resp.
canonical) determinant. One of these families of bundles is related
to the gradient of the $2\theta$ functions along the surface $C-C$,
the other family is the Brill-Noether locus $\mathcal{W}(3)$ of
stable bundles having at least three independent sections. We also
need (section 3) some relations between vectors of second order
theta functions, which one derives from Fay's trisecant formula
and its generalizations \cite{gu1}.
Finally (section 2.5) we describe the relationship between
these bundles and the objects discussed in \cite{geem-iz}, which are related
to the embedded tangent space at the origin to the moduli of stable
rank 2 bundles.

\bigskip

A natural question to ask is whether the equivalence \eqref{frobenius}
admits generalizations. 
We denote by $C_2$ the second symmetric
product of the curve. 
We introduce the fourfold $\cmc$, 
defined to be the image of the difference map $\phi_2 : \ctc
\ra JC$, which maps a $4$-tuple $(p+q,r+s)$ to the line bundle
$\mathcal{O}(p+q-r-s)$. 
With this notation we define the subseries
\begin{align}
\ggf & =\{ D \in \gg \: | \: \cmc \subset D \} \\
\ggg & =\{ D \in \gg \: | \: \mult_0(D) \geq 6 \}
\end{align}
and we prove the following theorems (section 5 and 6). Let $I(2)$ (resp.
$I(4)$) be the space of quadrics (resp. quartics) in canonical space
$|K|^*$ containing the canonical curve.

\begin{thm}
For any {non-trigonal} curve, there exists a canonical
isomorphism
\begin{equation} \label{exactsequence}
\ggd/\ggf \map{\sim} \sym^2 I(2).
\end{equation}
\end{thm}

\begin{thm}
For any non-hyperelliptic curve, we have the following inclusions
\begin{equation}\label{inclusions}
\ggf \subset \ggg \subset \ggd .
\end{equation}
The quotient of the first two spaces is isomorphic to the space of
quadratic syzygies among quadrics in $I(2)$, i.e.
\begin{equation} \label{quadsyzy}
\ggg / \ggf \cong \ker m : \sym^2 I(2) \lra I(4)
\end{equation}
where $m$ is the multiplication map.
\end{thm}

The deepest statement in Theorem 1.2 is the surjectivity  of the map 
in \eqref{exactsequence}. The proof uses essentially two ideas: first,
we can give an explicit basis of quadrics in $I(2)$ of rank less than
or equal to $6$ (Petri's quadrics, section 5.1) and secondly, we can
construct out of such a quadric a rank $2$ vector bundle having at least
four independent sections. This construction \cite{brivio} is recalled in 
section 5.2. Thus the Brill-Noether locus $\mathcal{W}(4)$ appears in
a natural way.

The proof of Theorem 1.3 is more in the spirit of Gunning's previous work and
uses only linear relations between vectors of second order theta functions.
Except for a few cases (section 6.2) we are unable to deduce the dimension
of $\ggg$ from \eqref{quadsyzy}.

In section 7 we give the version of Theorem 1.2 for trigonal curves.

It turns out that the vector bundle constructions used in the proofs
of Theorems 1.1 and 1.2 can be seen as examples of a global construction 
(section 8).

\bigskip

{\em Acknowledgements.} The authors are grateful to J. Fay, W.M. Oxbury and
S. Verra for various helpful comments. We especially thank B. van Geemen
for many fruitful discussions and his constant interest in this work.

\bigskip

{\em Notation.}

\begin{itemize}
\item If $X$ is a vector space or a vector bundle, by $X^*$ 
we denote its dual.
\item $K$ is the canonical divisor of the curve $C$.
\item For a vector bundle $E$ over $C$, 
$H^i(C,E)$ is often abbreviated by $H^i(E)$ and
$h^i(E) = \dim H^i(C,E)$.
\item $C_n$ is the $n$-th symmetric product of the curve $C$.
\item $W_d^r(C)$ is the subvariety of $Pic^d(C)$ consisting of line bundles
$L$ such that $h^0(L) > r$.
\item The canonical curve $C_{can}$ is the image of the embedding 
$\varphi_K : C \lra |K|^*$.
\item The vector space $I(n)$ is the space of degree $n$ hypersurfaces in
  $|K|^*$ containing $C_{can}$.
\end{itemize}

\section{Rank 2 vector bundles}

In this section we construct two families of stable rank 2
bundles over $C$ and describe some of their properties.

\subsection{Preliminaries on extension spaces}

Let $\mod0$ (resp. $\modk$) be the moduli space of rank 2 vector
bundles over $C$ with fixed trivial (resp. canonical) determinant. It can be 
shown
\cite{beau} that the Kummer maps given by the linear system
$|2\Theta_0|$ over $JC$ (resp. $|2\Theta|$ over $\pic^{g-1}(C)$) can
be factorized through the moduli space $\mod0$ (resp.$\modk$). This gives
the two following commutative diagrams
\begin{equation} \label{defD}
\begin{CD}
JC @>Kum>> |2\Theta_0|^{\ast} \\
@VViV            @VVwV \\
\mod0 @>D>> |2\Theta|
\end{CD}
\qquad
\begin{CD}
\pic^{g-1}(C) @>Kum>> |2\Theta|^{\ast} \\
@VViV            @VVwV \\
\modk @>D>> |2\Theta_0|
\end{CD}
\end{equation}
The vertical morphisms $i$ map $JC$ (resp. $\pic^{g-1}$) to the
semi-stable boundary of the moduli space $\mod0$ (resp. $\modk$)
by sending a line bundle $\xi$ to the split bundle $\xi \oplus \xi^{-1}$
(resp. $\xi \oplus K\xi^{-1}$). The rightmost morphism $D$ associates
to a semi-stable rank 2 bundle $E$ with canonical determinant a divisor $D(E)$,
whose support is given by
\begin{equation} \label{defd}
D(E) = \{ \xi \in JC  \: | \: h^0(E \otimes \xi) > 0 \}
\end{equation}
The definition of $D(E)$ for $E$ with trivial determinant is obtained
by substituing $JC$ by $\pic^{g-1}(C)$. The composite map $w \circ Kum =
D \circ i$ in the rightmost diagram \eqref{defD} is the translation morphism
\begin{eqnarray*}
\iota : \pic^{g-1}(C) & \lra & |2\Theta_0| \\
 \xi & \lms & \Theta_\xi + \Theta_{K\xi^{-1}}
\end{eqnarray*}
where the divisor $\Theta_\xi \subset JC$ is obtained by translating
the Riemann theta divisor $\Theta$ by $\xi$. Dually, we get $\iota: JC \lra
|2\Theta|$ by substituing $\Theta$ by $\Theta_0$. Note that the 
symmetric theta divisor $\Theta_0$ depends on the choice of a 
theta characteristic.

Given a subspace $V \subset H^0(\pic^{g-1}(C),2\Theta)$, we denote by
$V^\perp \subset H^0(JC,2\Theta_0)$ the image under the Wirtinger
duality \eqref{wirtinger} of its annihilator.

\subsubsection{degree 1}

Given an $x \in \pic^1(C)$, let $\mathbb{P}(x) = |Kx^2|^{\ast} =
\mathbb{P}H^1(C,x^{-2})$. This $g$-dimensional projective space
parametrizes isomorphism classes of extensions
\begin{equation} \label{llll}
 0 \lra x^{-1} \lra E \lra x \lra 0
\end{equation}
and the classifying map $\psi : \mathbb{P}(x) \ra \mod0 \ra |2\Theta|$
is linear and injective (lemme 3.6 \cite{beau2}). 
The following lemma (Prop. 1.2 \cite{prev})
describes the incidence relations between the extension
spaces $\mathbb{P}(x)$.
\begin{lem}
Let $x,y \in \pic^1(C)$. If $x \otimes y = \mathcal{O}(p+q)$, the intersection
$\mathbb{P}(x) \cap \mathbb{P}(y)$ is the secant line $\overline{pq}$
of the curve in either $\mathbb{P}(x)$ or $\mathbb{P}(y)$.
\end{lem}
The linear system $|Kx^2|$ can be used to map the curve into
$\mathbb{P}(x)$
$$\varphi = \varphi_{Kx^2} : C \longrightarrow \mathbb{P}(x).$$
A point in $\mathbb{P}(x)$ represents a stable bundle precisely away from
$C_x := \varphi(C)$, while the image of a point $p \in C $
represents the equivalence class of the semi-stable bundle
$x(-p) \oplus x^{-1}(p)$. The Abel-Jacobi map $t_x : C \ra
JC$ defined by $p \mapsto x(-p)$ fits in the following commutative
diagram
\begin{equation*}
\begin{CD}
C @>{t_x}>> JC \\
@VV{\varphi}V   @VViV \\
\mathbb{P}(x) @>{\psi}>> \mod0
\end{CD}
\end{equation*}

\subsubsection{degree 2}

Given an $x \in \pic^2(C)$, we define similarly $\mathbb{P}(x) =|Kx^2|^{\ast}
= \mathbb{P} H^1(C,x^{-2})$. This $(g+2)$-dimensional projective
space parametrizes
as in section 2.1.1 isomorphism classes of extensions \eqref{llll}
and the classifying map  $\psi : \mathbb{P}(x) \ra \mod0 \ra |2\Theta|$
is given by the complete system of quadrics \cite{bertram} through the 
embedded curve
$\varphi_{Kx^2}: C \hookrightarrow \mathbb{P}(x)$. Thus this map is defined 
away from
the curve $C$. We will use the fact that a secant line
$\overline{pq}$ to the curve $C \subset \mathbb{P}(x)$ is contracted by
the linear system of quadrics to a point represented by the split bundle
$x(-p-q) \oplus x^{-1}(p+q)$.

\subsection{The bundles $E(p,q,r)$}

We will associate to any three distinct points $p,q,r$ on the
curve $C$ a unique stable rank 2 bundle $E(p,q,r)$, which we identify
with its image in $|2\Theta|$. Consider the line bundle
\begin{equation} \label{definitionx}
x =  {\cal O}(p + q -r) \in \pic^1(C)
\end{equation}
and denote by $T_{p}(C_x)$ the embedded tangent line to the curve $C_x \subset
\mathbb{P}(x)$
at the point $p \in C$. Note that the image $C_x \subset \mathbb{P}(x)$
is smooth, since $h^0(x^2) = 0 $ (Lemma 5.2 \cite{lange}).
Recall that $\mathbb{P}(\mathbb{T}_0) \subset | 2 \Theta |$ denotes the
linear span of the surface $C-C \subset JC$.

We say that three points $p,q,r$ are collinear if they are collinear
as points  on the canonical curve $C_{can}$.

\begin{lem}
$T_p(C_x) \subset \mathbb{P}(\mathbb{T}_0) \; \iff \; p,q,r$ are collinear
\end{lem}
\begin{proof}
Note that the image $t_x(p)  = {\cal O}(q-r) \in C-C
 \subset JC$. We identify the projectivized
tangent space to $JC$ at the point $\varphi(p)$ with the canonical
space  $|K|^{\ast}$. Under this identification the projectivized
tangent space to $C-C$ at $q-r$  corresponds to the
secant line $\overline{qr} \subset |K|^{\ast}$
and the projectivized tangent line to the curve
$t_x(C) \subset JC$ at the point $q-r$ corresponds to
$p \in |K|^{\ast}$. If $p \in \overline{qr}$,
then $T_p(C_x) \subset T_{q-r}(C-C) \subset \pp(\mathbb{T}_0)$.

Conversely, by \cite{izadi2}
Th\'eor\`eme A, 
 we know that
the projectivized tangent spaces to divisors in $\gg
= \mathbb{T}_0^{\bot}$ at the point $q-r$ cut out the secant line
$\overline{qr} \subset |K|^{\ast}$. Therefore,
$T_p(C_x) \subset \pp(\mathbb{T}_0)$ implies that
$p \in \overline{qr}$.
\end{proof}
\bigskip
We now introduce the line bundles
\begin{equation} \label{definitionyz}
y = {\cal O}(p + r - q) \hspace{2cm} z= {\cal O}(q+ r - p)
\end{equation}

\begin{prop}
For any non-collinear points $p,q,r$, there exists a unique
stable bundle $E := E(p,q,r) \in \mod0$ containing the
three line bundles $x^{-1}, y^{-1}, z^{-1}$. Moreover,
$E \notin \mathbb{P}(\mathbb{T}_0)$.
\end{prop}

\begin{proof}
By \cite{lange} Prop. 1.1, the extension classes
lying on the tangent line $T_p(C_x) \subset \mathbb{P}(x)$
correspond to stable rank $2$ bundles containing $x(-2p)
= y^{-1}$. By Riemann-Roch
$$h^0(Kx^2(-2p -2q)) = h^0(K(-2r)) = g-2$$
\noindent
hence the two tangent
lines $T_p(C_x)$ and $T_q(C_x)$ are contained
in a projective plane and are distinct. Their
intersection point $T_p(C_x) \cap T_q(C_x)$ determines a stable bundle, 
which we
denote by $E(p,q,r)$ and which contains, by
construction, $x^{-1}, y^{-1}$ and $z^{-1}$. We now observe that
$$x \otimes y = {\cal O}(2p)$$
\noindent
hence, by Lemma 2.1, $T_p(C_x) = T_p(C_y)$.
Similarly $T_q(C_x) = T_q(C_z)$ and $T_r(C_z) =
T_r(C_y)$, from which it is clear that the three
tangent lines intersect in a single point, hence the
construction does not depend on the order of the three points.
Finally Lemma 2.2 implies that $E \notin \mathbb{P}(\mathbb{T}_0)$.
\end{proof}

\begin{prop}
If $a \in \sing$, then $h^0(E(p,q,r) \otimes a) >0.$
\end{prop}

\begin{proof}
We write $E = E(p,q,r)$ and take $x$ as in \eqref{definitionx}.
We tensor the exact sequence \eqref{llll}
with $a$ and take global sections
$$0 \lra H^0(x^{-1} \otimes a) \lra H^0(E \otimes a) \lra H^0(x \otimes a)
\map{\delta} H^1(x^{-1} \otimes a)$$
If $h^0(x^{-1} \otimes a) >0$, we are done. So we can assume that
$h^0(x^{-1} \otimes a) = 0$. By Serre duality, we have $h^0(E \otimes a) 
= h^0(E \otimes a')$ with $a' = K \otimes a^{-1}$, so we can also
assume that $h^0(x^{-1} \otimes a') = 0$. But this is equivalent
to $h^0(x \otimes a) = h^1(x^{-1} \otimes a) = 1$. Finally we see that
$h^0(E \otimes a) >0$ if and only if the coboundary map $\delta$ is
zero, i.e. the one-dimensional image of the multiplication map
$$H^0(x \otimes a) \otimes H^0(x \otimes a') \lra H^0(Kx^2)$$
is contained in the hyperplane of $H^0(Kx^2)$ corresponding to
the bundle $E$ and which is, by definition, the linear span of
the two subspaces $H^0(Kx^2(-2p))$ and $H^0(Kx^2(-2q))$. We observe
now that $H^0(a(-r)) \subset H^0(x \otimes a)$ and for dimensional
reasons these two spaces must coincide, hence the unique global
section of $x \otimes a$ vanishes  at $p$ and $q$. Since the
same holds for $x \otimes a'$, we are done.
\end{proof}

\subsection{The bundles $E_W$}

In this section we recall some constructions from \cite{geem-iz}.
Let $Gr(3,H^0(K))$ be the Grassmannian of 3-planes in
$H^0(K)$ and $\mathcal{W}(3)$ the locus of stable bundles $E$ in $\modk$ 
that are generated by global sections
and for which $h^0(E) = 3$.
We will define
a rational map

$$Gr(3,H^0(K)) \map{\beta} \mathcal{W}(3). $$
\noindent
For a generic 3-plane $W\subset H^0(K)$, it can be shown that the 
multiplication map
\begin{equation}  \label{mult}
W \otimes H^0(K) \rightarrow H^0(K^2) \ \ \ {\rm is\ \  surjective.}
\end{equation}
\noindent
Then the dual $E_W^{\ast}$ of the bundle $\beta(W) = E_W$ is
defined by the exact sequence
\begin{equation} \label{definitionew}
0 \longrightarrow E_W^{\ast} \longrightarrow W \otimes
{\cal O}_C \map{ev} K \longrightarrow 0
\end{equation}
\noindent
where the last map is the evaluation map of global sections. Then
condition \eqref{mult} implies that $E_W \in \mathcal{W}(3)$.

Conversely, we can associate to a bundle
$E \in \mathcal{W}(3)$ a 3-plane $W_E$, i.e. we have an {\it inverse} map
$$\mathcal{W}(3) \map{\alpha} Gr(3,H^0(K)).$$
\noindent
To define $\alpha(E) = W_E$, we consider the exact sequence
$$ 0 \longrightarrow K^{-1} \map{i} H^0(E) \otimes \mathcal{O}_C
\map{ev} E \longrightarrow 0$$
\noindent
then the dualized exact sequence induces an injective
linear map on global sections
$$ H^0(i^{\ast}) \ : \ H^0(E)^{\ast} \longrightarrow H^0(K)$$
\noindent
and we let $W_E = {\rm im} \ H^0(i^{\ast})$.

\noindent
Moreover, under the natural duality $H^0(E)^{\ast} \cong \bigwedge^2 H^0(E)$,
$W_E$ coincides with the image of the exterior product map
$\bigwedge^2 H^0(E) \rightarrow H^0(K)$.

\begin{rem}
Since $E_W$ is generated by global sections we can define a map
\begin{eqnarray*}
C & \lra & \mathbb{P}H^0(E_W) = \mathbb{P}^2 \\
p & \lms & s_p := \ker (H^0(E_W) \xrightarrow{ev} {E_W}_{|p})
\end{eqnarray*}
which associates to a point $p$ the unique section $s_p$ of $E_W$ vanishing
at $p$. This map coincides with the canonical map $\varphi_K : C
\hookrightarrow |K|^{\ast}$ followed by the projection with center
$\mathbb{P}W^{\bot} \subset |K|^{\ast}$.
\end{rem}

\subsection{Incidence relations}

In the previous sections we have constructed two families of bundles
i.e. $E(p,q,r)$ (section 2.2) and $E_W$ (section 2.3); each gives
$2\theta$-divisors under the $D$-maps in $|2\Theta|$ and $|2\Theta_0|$.
These spaces are dual to each other \eqref{wirtinger} and it is
therefore useful to determine their incidence relations. We will denote by
$H_W \in |2\Theta|^{\ast}$ the hyperplane in $|2\Theta|$ corresponding
under \eqref{wirtinger} to $D(E_W) \in |2\Theta_0|$ and by $p \wedge q
\wedge r \in \Lambda^{3} H^0(K)^{\ast}$ (resp. $\Lambda^3 W \in \Lambda^3 
H^0(K)$) \ the Pl\"ucker image of the
3-plane in $H^0(K)^{\ast}$ spanned by the 3 points $p,q,r \in C_{can}$ 
(resp. the 3-plane $W$ in $H^0(K)$).

\begin{prop}
We have the following equivalence
$$E(p,q,r) \in H_W \: \: \iff \: \: \langle \Lambda^3 W ,
p \wedge q \wedge r \rangle = 0. $$
\end{prop}

\begin{proof}
We write $E= E(p,q,r)$ and take $x$ as in \eqref{definitionx}. We
note that the first condition $E \in H_W$ is equivalent to
$h^0(E_W \otimes E) > 0$ \cite{beau2}. Using remark 2.5, we see that the
second condition is equivalent to the three sections $s_p,s_q,s_r \in
\mathbb{P}H^0(E_W) = \mathbb{P}^2$ being collinear.

We tensor the exact sequence \eqref{llll} with $E_W$ and take global
sections
\begin{equation} \label{exasqew}
0 \lra H^0(E_W \otimes x^{-1})  \lra H^0(E_W \otimes E) \lra 
 H^0(E_W \otimes x) \map{\delta(\epsilon)} H^1(E_W \otimes x^{-1})
\end{equation}

Let us first consider the case when $h^0(E_W \otimes x^{-1}) >0$. Then
$h^0(E_W \otimes E) > 0$, and, by definition of the
bundle $E_W$ \eqref{definitionew}, we have
\begin{equation} \label{xxx}
H^0(E_W \otimes x^{-1}) = \ker (W \otimes H^0(Kx^{-1}) \map{ev}
H^0(K^2x^{-1}))
\end{equation}
But $H^0(Kx) = H^0(K(-p-q))$. Furthermore, if $\dim \: W \cap
H^0(K(-p-q)) = 1$, then the kernel of the map \eqref{xxx} is zero.
Hence $\dim \: W \cap H^0(K(-p-q)) \geq 2$, which implies that
$\dim \: W \cap H^0(K(-p-q-r)) \geq 1$ and we are done.
\bigskip

Now we can assume that $h^0(E_W \otimes x^{-1})= 0$, or equivalently
$h^0(E_W \otimes x) = 2$. Next we observe that the coboundary map
$\delta(\epsilon)$ \eqref{exasqew}, which depends on the extension
class $\epsilon \in \mathbb{P}(x)$ of the bundle $E$, is skew-symmetric
(here we identify $H^1(E_W \otimes x^{-1}) = H^0(E_W \otimes x)^{\ast}$).
Therefore, if $\{s,t\}$ is a basis of $H^0(E_W \otimes x)$, we have
\begin{equation} \label{xxxx}
\begin{split}
 h^0(E_W \otimes E) >0 \: \: & \iff  \: \: \det (\delta(\epsilon)) = 0 \\
 & \iff  \: \: s \wedge t \in H^0(Kx^2) \: \text{vanishes doubly at $p$ and 
$q$.}
\end{split}
\end{equation}
Let $\sigma$ be a section of the line bundle $\mathcal{O}(p+q)$. The
4-dimensional space $H^0(E_W(p+q))$ contains $H^0(E_W)$ and
$H^0(E_W \otimes x)$, which intersect in the section $s_r \otimes \sigma$.
We can choose $s = s_r \otimes \sigma$ and $t \notin H^0(E_W)$.
Then $t\wedge s_p$ vanishes at $p$ and $r$, hence $t \wedge s_p
\in H^0(K(q-r)) = H^0(K(-r))$, so $t \wedge s_p$ also vanishes at $q$.
Similarly, $t \wedge s_q$ vanishes at $p$. The condition \eqref{xxxx}
says that $t \wedge s_r$ vanishes at $p$ and $q$. Since $t \notin
H^0(E_W)$ we can assume e.g. that $t(q) \not= 0$, but then
$s_p(q) \wedge s_r(q) = 0$, which implies that the three sections
$s_p,s_q,s_r$ cannot be linearly independent (otherwise they would
generate $E_W$ at $q$).
\end{proof}

\begin{rem}
The incidence relations were first proved in \cite{geem-iz}
for slightly different objects (see section 2.5). Working with the
bundles $E(p,q,r)$ instead of the projective spaces $\pp^4_{p,q,r}$
simplifies the proof somewhat.
\end{rem}

\subsection{Other descriptions}

Consider the projection $P$ with center $\mathbb{P}(\tankum)$
\begin{equation} \label{projecti}
 P  : |2 \Theta| \lra \mathbb{P}(H^0(2\Theta)/\tankum)
= |2 \Theta|_{\tankum}
\end{equation}
In section 2.2 we have associated to any triple of
non-collinear points $p,q,r$
a point in $|2 \Theta|_{\tankum}$, namely $P \circ D(E(p,q,r))$.
In \cite{geem-iz} the authors associate to the same data a
4-dimensional projective space $\mathbb{P}^4_{p,q,r} \subset |2\Theta|$.
Their construction goes as follows:

We take $\zeta = \mathcal{O}(p+q) \in \pic^2(C)$ and consider in
$\mathbb{P}(\zeta)$ (see 2.1.2) the 3-dimensional subspace $\langle 2p+2q+r 
\rangle$
spanned by the tangent lines at $p$ and $q$ to $C \hookrightarrow
\mathbb{P}(\zeta)$ and the point $r$. Then the restricted linear
system of quadrics through $C$ determines a rational map
$$\psi : \langle 2p+2q+r \rangle \lra |2\Theta|.$$
The image of $\psi$ is a cubic threefold in $\mathbb{P}^4_{p,q,r}
\subset |2\Theta|$, singular at $D(\mathcal{O} \oplus \mathcal{O})$.
Furthermore if $p,q,r$ are non-collinear, $\dim \pp(\tankum)
\cap \mathbb{P}^4_{p,q,r} = 3$, hence the projective space
$\mathbb{P}^4_{p,q,r}$ is contracted by the projection $P$ to a point.

\begin{lem}
For any non-collinear points $p,q,r$
$$P \circ D(E(p,q,r)) = P (\mathbb{P}^4_{p,q,r})
\in |2 \Theta|_{\tankum}$$
\end{lem}

\begin{proof}
We need to show that there exists an extension class $\epsilon \in
\langle 2p+2q+r \rangle \subset \mathbb{P}(\zeta)$ which corresponds to the
bundle $E(p,q,r)$ mod $\mathbb{P}({\tankum})$. We recall that
$E(p,q,r)$ has been characterized in Prop. 2.3. and we note
(\cite{lange} Prop. 1.1) that a bundle $E(\epsilon)$ fitting in the sequence
$$ 0 \lra \zeta^{-1} \lra E(\epsilon) \lra \zeta \lra 0$$
contains $\zeta(-2q-r) = z^{-1}$ if and only if 
$\epsilon \in \langle 2p+2q+r \rangle$.
Similarly $E(\epsilon)$ contains $\zeta(-2p-r)= y^{-1}$ if and only if 
$\epsilon
\in \langle 2p+r \rangle$. Therefore we may take $\epsilon \in \langle 2q+r 
\rangle \cap
 \langle 2p+r \rangle \not= \emptyset$, so that $E(\epsilon)$ contains $y^{-1}$
and $z^{-1}$. Such extension classes $\epsilon$ are parametrized by
$\mathbb{P}(y) \cap \mathbb{P}(z)$, which is the tangent line
$T_r(C_z) = T_r(C_y)$ at the point $\mathcal{O}(p-q) \in
\mathbb{P}(\tankum)$, hence $P(E(\epsilon)) =
P(E(p,q,r))$ for any such $\epsilon$, and we are done.

\end{proof}

\section{Gunning's results on second order theta functions}

In this section we recall some classical theory of theta functions 
seen as holomorphic quasi-periodic
functions on $\mathbb{C}^g$, as well as some results by Gunning on the
gradient and Hessian of $2\theta$-functions along the surface $C-C$. 
We refer to Fay's book \cite{fay} and to
\cite{gu1} for a detailed exposition.

Let $\ca$ be the universal covering space of the curve
$C$. We choose a base point $z_0 \in \ca$ and a canonical set
of generators of $H_1(C,\mathbb{Z})$ and call the
corresponding canonical basis of Abelian differentials
$\omega_1,\ldots,\omega_g \in H^0(C,K)$; these can be
thought of as holomorphic differential 1-forms on $\ca$
invariant under the group $\Gamma$ of covering tranformations
acting on $\ca$. We construct from
these data the period matrix $\Omega$. The associated
Abelian integrals
$$ w_j(z) = \int_{z_0}^{z} \omega_j $$
\noindent
are holomorphic functions on $\ca$ and are the coordinate
functions for a map $w : \ca \lra \mathbb{C}^g$.
Let $\Lambda$ be the lattice in $\mathbb{C}^g$ defined by the period
matrix $\Omega$, then $JC = \mathbb{C}^g/\Lambda$ and we have a
commutative diagram
$$
\begin{CD}
\ca @>w>>  \mathbb{C}^g \\
@VV{\pi}V   @VVV \\
\ca/\Gamma = C @>t_x>> JC
\end{CD}
$$
The horizontal map $t_x$ is the Abel-Jacobi map (see 2.1.1) with
$x= \mathcal{O}(p)$ and $p = \pi(z_0)$.
Both vertical arrows are quotient maps of the group actions
of $\Gamma$ (acting on $\ca$) and $\Lambda$ (acting on $\cc^g$). 
Sections of the
line bundle $\mathcal{O}(2 \Theta_0)$ over $JC$ correspond to the
classical second-order theta functions. A basis of the space
$H^0(JC, \mathcal{O}(2 \Theta_0))$ is given by the
holomorphic functions on $\cc^g$
\begin{equation} \label{twotheta}
\theta_2 \thetachar{\nu}{0} (w,\Omega) = \theta \thetachar{\nu}{0} (2w,2\Omega)
\  \ {\rm for} \ \nu \in (\frac{1}{2} \mathbb{Z} / \mathbb{Z})^g
\end{equation}
\noindent
where the right-hand side is obtained from the first-order
theta function with characteristics  $\thetachar{\nu}{0}$
$$
\theta \thetachar{\nu}{0} (w,\Omega) = \sum_{n \in \mathbb{Z}^g}
\exp 2\pi i \left[ \frac{1}{2} {}^t(n+\nu)\Omega(n+\nu) +
(n+\nu)w \right]
$$
\noindent
The functions \eqref{twotheta} are the coordinate
functions of a holomorphic map
\begin{eqnarray} \label{kummermap}
\begin{array}{rcl}
\theta_2 : \mathbb{C}^g & \longrightarrow & \mathbb{C}^{2^g} \\
 w & \longmapsto & (\ldots ,\theta_2 \thetachar{\nu}{0}(w,\Omega),\ldots)
\end{array}
\end{eqnarray}
\noindent
Using this basis, we identify  (via the Wirtinger duality \eqref{wirtinger})
$\mathbb{P}( \mathbb{C}^{2^g}) = |2 \Theta|$, so that the map
\eqref{kummermap} coincides
with the Kummer map $JC  \rightarrow |2  \Theta|$.

\bigskip
Next we introduce the {\em prime form} (see  \cite{gu1} formula (22))
$q(z_1,z_2)$, which is a holomorphic function on $\ca
\times \ca$ with a simple zero along the subvariety
$z_1  = T z_2$ for all covering transformations $T$ and vanishing
nowhere else. Moreover, $q(z_1,z_2) = -q(z_2,z_1)$.  Alternatively
the function $q$ is, up to a multiplicative constant, the pull-back 
to $\ca \times \ca$
of the unique global section having as zero scheme  the
diagonal in $C \times C$.

We use canonical coordinates on the universal covering $\ca$ (see section 6
\cite{gu1}) and we denote by $w'_j$ the derivative of the holomorphic function
$w_j$ with respect to the canonical coordinates. To a point $a \in \ca$ we 
associate
the differential operator 
$$ D_a = \sum_{j=1}^g w'_j(a) \frac{\partial \:}{\partial w_j}.$$
This operator corresponds, up to multiplication by a scalar, to the unique
translation-invariant vector field over $JC$, which has as tangent vector at 
the 
origin $O \in JC$ the tangent direction at $O$ to the curve $C$, where
$C$ is embedded in $JC$ by $q \mapsto \mathcal{O}(q-p)$ with $p = \pi (a) \in 
C$.

Gunning (see \cite{gu1} formulae (41), (42), (44)) introduces for
$a_1,a_2,a_3,a_4 \in \ca$ the following vectors in $H^0(2\Theta)/\tankum = 
\cc^{2^g}/\tankum$
\begin{equation}
\xi(a_1,a_2,a_3)  =  q(a_2,a_3)^{-2} P D_{a_1}\theta_2(w(a_2 - a_3))
\end{equation}
\begin{equation} \label{defsigma}
\sigma(a_1;a_2,a_3;a_4)  =  \left[ \frac{\partial \:}{\partial a_1} \log
\frac{q(a_1,a_2)}{q(a_1,a_3)} \right] \cdot \xi(a_2,a_3,a_4) 
\end{equation}
\begin{equation}
\tau(a_1,a_2;a_3,a_4)  =  q(a_3,a_4)^{-2} P D_{a_1} D_{a_2} \theta_2(w(a_3 - 
a_4))
\end{equation}
Then $\xi$ (resp. $\tau$) defines a holomorphic function on $\ca^3$ (resp. 
$\ca^4$)
with values in $\cc^{2^g}/\tankum$ and $\sigma$ defines a meromorphic function
on $\ca^4$ which has as singularities at most simple poles along the loci
$a_1 = Ta_2$ and $a_1 = Ta_3$ for all covering transformations $T \in \Gamma$.
The functions $\xi,\sigma,\tau$ have the following symmetry properties:
$\xi$ is skew-symmetric in $a_1,a_2,a_3$ (Cor. 4 \cite{gu1}), $\sigma$ is 
symmetric
in $a_2,a_3$ and $\tau$ is symmetric in $a_1,a_2$ and in $a_3,a_4$.
Given four points $a_1,a_2,a_3,a_4 \in \ca$ we will let $\sigma_{ijkl} = 
\sigma(a_i;a_j,a_k;a_l)$, $\tau_{ijkl} = \tau(a_i,a_j;a_k,a_l)$ and
$q_{ij} = q(a_i,a_j)$  where $i,j,k,l$
are four indices such that $\{i,j,k,l\} = \{1,2,3,4\}$. Then 
 the following statements hold:
\begin{thm}[Theorem 3 \cite{gu1}]
There exist vectors $\xi_{jkl} \in \cc^{2^g}/{\tankum}$ which are
skew-sym- metric in their indices j,k,l and such that
$$ \xi(a_1,a_2,a_3) = \sum_{j,k,l} \xi_{jkl} w'_j(a_1) w'_k(a_2)
w'_l(a_3). $$
\end{thm}
\begin{cor}
The dimension of the linear span of the vectors
$\xi(a_1,a_2,a_3)$ for $a_i$ varying in $\ca$ is at most 
$\binom{g}{3}$.
\end{cor}
\begin{prop}[Cor. 5 \cite{gu1}]
For any points $a_1,a_2,a_3,a_4 \in \ca$
\begin{enumerate}
\item[(1)]
 $\frac{1}{2} PD_{a_1}D_{a_2}D_{a_3}D_{a_4}\theta_2(0) = \tau_{1234}
+ \tau_{1324} + \tau_{2314} - 2\sigma_{1234} - 2\sigma_{2134} -2\sigma_{3124}$
\item[(2)]
$\frac{1}{2} \tau_{1324} +\frac{1}{2} \tau_{1423} = \left( \frac{q_{12}q_{34}}
{q_{13}q_{14}q_{23}q_{24}} \right)^2 P\theta_2(w(a_1+a_2-a_3-a_4))
+ \sigma_{1342} - \sigma_{3241} - \sigma_{4231}$
\end{enumerate}
\end{prop}

\section{Proof of theorem 1.1}

We consider the rational map
\begin{eqnarray*}
\rho : C_3 & \lra & |2\Theta|_{\tankum} \\
(p,q,r) & \longmapsto & P \circ D(E(p,q,r))
\end{eqnarray*}
By Lemma 2.2 the rational map $\rho$ is defined away from the 
triples of collinear points. In particular, $\rho$ is a morphism
if $C$ is not trigonal. 
Let $\mathbb{P}(\tanmod)$ be the inverse image under the projection $P$
\eqref{projecti}  of the linear
span of the image of $\rho$. Obviously $\tankum \subset \tanmod$.
Since the bundle $E(p,q,r) \in T_p(C_x)$, the tangent line at the
point $\mathcal{O}(q-r) \in \mathbb{P}(\tankum)$ to the curve
$C_x$, we see  that a hyperplane in $|2 \Theta|$ containing
$\pp(\tankum)$ and the point $D(E(p,q,r))$ determines via \eqref{wirtinger}
a $2\theta$-function $f$ such that $D_p f(q-r) = 0$. Hence, varying the
points $p,q,r$, we obtain the
equality $\tanmod^{\perp} = \ggd$.

\bigskip

The rest of the argument coincides with the argument given in
 \cite{geem-iz}. We therefore just sketch their proof: the rational map
$\rho$ factorizes as follows
\begin{eqnarray}
\begin{array}{ccc}
C_3 & \map{\rho} & \mathbb{P}(\tanmod/\tankum)\\
 & & \\
\downarrow^{\pi} & & \uparrow^{\gamma} \\
 & & \\
Gr(3,H^0(K)^{\ast}) & \map{Pl} & \mathbb{P}(\Lambda^3 H^0(K)^{\ast})
\end{array}
\end{eqnarray}
where $\pi(p,q,r) = p \wedge q \wedge r$ and $Pl$ is the Pl\"ucker embedding.
By the incidence relations (Prop. 2.6) the support of the pull-back 
$\rho^{\ast} H_W$
equals the set
$$ \{ (p,q,r) \in C_3 \: | \: \langle \Lambda^3 W,  p \wedge q
\wedge r \rangle = 0  \}. $$
More precisely, we have an equality $\rho^{\ast}
H_W = (Pl \circ \pi)^{\ast}((\Lambda^3 W)^{\perp})$ as divisors on
$C_3$, hence the factorization and the injectivity of $\gamma$.
Since $\dim (\tanmod / \tankum) \leq \binom{g}{3}$ by Cor. 3.2,
$$\gamma : \Lambda^3 H^0(K)^* \map{\sim} \tanmod/\tankum$$
is an isomorphism. Then we take duals and since 
$\tanmod^{\perp} = \ggd$ and $\tankum^{\perp} = \gg$ we obtain Thm. 1.1 (1).
\bigskip

To prove Thm. 1.1 (2), we note that we have an inclusion 
$\spansing \subset \ggd$, since by 
Riemann's singularity theorem, for $\xi \in \sing$
$$\mult_{p-q}(\Theta_\xi + \Theta_{K\xi^{-1}}) \geq 2h^0(\xi(p-q))
\geq 2.$$
Since both spaces have the same dimension (for a computation of $\dim 
\spansing$ 
see
\cite{geem-iz}), we obtain equality.

\section{Quadrics on canonical space}

\subsection{Petri's quadrics}
We will denote by $\tilde{Q}$ the polarized form of a quadric $Q$ on canonical
space $|K|^*$, i.e. $\tilde{Q}$ is the symmetric bilinear form such that
$\tilde{Q}(v,v) = Q(v)$,  $\forall v \in H^0(K)$.

\begin{lem}
We consider $g-2$ points in general position $p_1,\ldots,p_{g-2}$. If a 
quadric $Q \in I(2)$ is such that
$$ \tilde{Q}(p_i,p_j) = 0 \qquad \forall i,j \in \{ 1,\ldots,g-2 \} $$
then $Q$ is identically zero.
\end{lem}

\begin{proof}
By the general position theorem (\cite{acgh} p.109), we know that a general
hyperplane $H \subset |K|^*$ meets $C$ in $2g-2$ points any $g-1$ of which
are linearly independent. We consider such an $H$ and $g-1$ independent 
points $q_1,\ldots, q_{g-1}$ in $H$ and suppose that the $g-2$ points $p_1,
\ldots, p_{g-2}$ are taken among the $g-1$ residual points of $H \cap C$, i.e
$$ H \cap C = \{ p_1, \dots, p_{g-2}, q_1, \ldots, q_{g-1},q_g \} . $$
It is clear that a general $(g-2)$-tuple $(p_1,\ldots,p_{g-2})$ can be 
realized
in this way.

Now suppose that $\tilde{Q}(p_i,p_j) = 0$, $\forall i,j \in \{ 1, \ldots, g-2 
\}$,
i.e. $Q$ contains the linear subspace $\Pi$ spanned by the $p_i$'s. Since $\Pi$
is a hyperplane in $H$, any line $\overline{q_iq_j}$, for $1 \leq i < j 
\leq g$ intersects $\Pi$. So $\overline{q_iq_j}$ is entirely contained in $Q$,
since it meets $Q$ in at least three points. Hence we obtain that 
$\tilde{Q}(q_i,q_j) = 0$, $\forall i,j \in \{ 1, \ldots, g \}$, i.e. $Q$ 
contains the hyperplane $H$. But since $Q$ contains $C$, it cannot be the 
union of two hyperplanes, hence $Q$ is identically zero.
\end{proof}

\bigskip

 From now on we fix $g-2$ points $p_1, \ldots, p_{g-2}$, which are chosen in
general position. By the preceding Lemma 5.1, the $\binom{g-2}{2}$ hyperplanes
(for  $1 \leq i < j \leq g-2$)
$$ \mathcal{H}_{ij} = \{ Q \in I(2) \: | \: \: \tilde{Q}(p_i,p_j) = 0 \} $$
are linearly independent in $I(2)^*$, hence they form a basis of $I(2)^*$. Let 
us
denote by $\{ Q_{ij} \}$ the corresponding dual basis. The quadrics 
$Q_{ij}
\in I(2)$ are characterized by the properties
\begin{equation} \label{charquad}
\begin{split}
 & \tilde{Q}_{ij}(p_\alpha,p_\beta) = 0 \: \: \text{if} \: \{i,j\} \not= 
\{ \alpha,\beta \} \\
 & \tilde{Q}_{ij}(p_i,p_j) \not= 0.
\end{split}
\end{equation}

\bigskip

This basis of quadrics has been used by K. Petri (\cite{petri}, see also 
\cite{acgh}
p.123-135) in his work on the syzygies of the canonical curve. He defines them
in a slightly different way:

Choose two additional points $p_{g-1},p_g$ in general position. For each $i$,
$1 \leq i \leq g$, pick a generator $\omega_i$ of the one-dimensional space
\begin{equation} \label{defomega}
H^0(K(-\sum^g_{j=1,j \not= i} p_j)) = \mathbb{C} \omega_i 
\end{equation}

Up to a constant the $\omega_i$'s form a dual basis to the points $p_i \in 
|K|^*$.
Then there are constants (\cite{acgh} p.130) $\lambda_{sij}, \mu_{sij}, b_{ij} 
\in \mathbb{C}$ such that, if we let 
\begin{equation} \label{vect}
\eta_{ij} = \sum_{s=1}^{g-2} \lambda_{sij} \omega_s \qquad 
\nu_{ij} = \sum_{s=1}^{g-2} \mu_{sij} \omega_s
\end{equation}
\noindent
the quadratic polynomials
\begin{equation} \label{Petriquad}
R_{ij} = \omega_i \omega_j - \eta_{ij} \omega_{g-1} - \nu_{ij} \omega_g
- b_{ij} \omega_{g-1} \omega_{g}
\end{equation}
\noindent
all vanish on $C$. Moreover the $R_{ij}$'s form a basis of $I(2)$ and, 
obviously,
the rank of the quadric $R_{ij}$ is less than or equal to $6$.

\begin{lem}
For each $1 \leq i < j \leq g-2$, the quadrics $R_{ij}$ satisfy the conditions
\eqref{charquad}, hence $R_{ij} = Q_{ij}$.
\end{lem}
\begin{proof}
This follows immediately from the definition \eqref{defomega} of the 
$\omega_i$'s.
\end{proof}

\begin{prop}
If $C$ is neither trigonal nor a smooth plane quintic and the points $p_1, 
\ldots,
p_{g-2}$ are in general position, then the quadrics $Q_{ij}$ have the following
properties
\begin{align}
  (i) \:  &  \si Q_{ij} \cap C = \emptyset \\
  (ii) \: & \rk Q_{ij} = 5 \: \text{or} \: 6.
\end{align}
\end{prop}

\begin{proof}
We fix two indices $i,j$. First we observe that the singular locus $\si 
Q_{ij}$ is
the annihilator of the linear space 
\begin{equation} \label{annsing}
\langle \omega_i, \omega_j, \omega_{g-1}, \omega_{g}, \eta_{ij}, \nu_{ij} 
\rangle \subset |K|.
\end{equation}
Hence $C \cap \si Q_{ij}$ is the base locus of this linear subsystem. In 
particular
$C \cap \si Q_{ij}$ is contained in the base locus of $\langle \omega_i, 
\omega_j,
\omega_{g-1}, \omega_g \rangle$, which, by construction, consists of the $g-4$
points (we delete the $i$-th and $j$-th point)
\begin{equation} \label{points}
 p_1, \ldots, \hat{p}_i, \ldots, \hat{p}_j, \ldots, p_{g-2}.
\end{equation}
We will denote by $D_{ij}$ the degree $g-4$ divisor with support
\eqref{points} and by $\bar{D}_{ij}$ the linear span of $D_{ij}$ in
$|K|^*$.
Suppose now that there exists a $(g-2)$-tuple $p_1, \ldots, p_{g-2}$ such that
(i) holds, then, since (i) is an open condition, (i) holds for a general 
$(g-2)$-tuple of points. Therefore we will assume that (i) does not hold for
all $(g-2)$-tuples, i.e. $\forall p_1,\ldots,p_{g-2}$ (in general
position), there exists a $p_\alpha \in \si Q_{ij}$ for some $\alpha
\in \{ 1, \ldots,g-2\}, \alpha \not= i,j$. But, since the quadric $Q_{ij}$
does not depend on the order of the $g-4$ points \eqref{points}, $p_\alpha \in
\si Q_{ij}$ implies that all $g-4$ points \eqref{points} are in $\si Q_{ij}$.
Hence
\begin{equation} \label{pointssi}
\bar{D}_{ij} \subset \si Q_{ij}
\end{equation}
and therefore $\rk Q_{ij} \leq 4$. Since 
$\omega_i, \omega_j, \omega_{g-1}, \omega_g$ are linearly independent, we
have $\rk Q_{ij} =4$. Hence the inclusion \eqref{pointssi} is an equality (same
dimension).

Consider now the two rulings of the rank $4$ quadric $Q_{ij}$: they cut out on 
the 
curve two pencils of divisors
$$ \mathbb{P}^1 \subset |L| \qquad \mathbb{P}^1 \subset |M| $$
\noindent
such that $L,M$ are line bundles satisfying 
\begin{equation} \label{relationlb}
L \otimes M = K(-D_{ij}).
\end{equation}
Therefore for general $D_{ij}$, we have constructed a pair $(L,M) \in 
W^1_d(C) \times W^1_{d'}(C)$ satisfying \eqref{relationlb} with $d + d' = 
\deg K(-D_{ij}) = g+2$. By Mumford's refinement of Martens' Theorem (see 
\cite{acgh}
p.192-3), if $C$ is neither trigonal, bi-elliptic, nor a smooth plane quintic,
then $\dim W^1_d(C) \leq d-4$ for $4\leq d \leq g-2$. Hence
$$\dim W^1_d(C) \times W^1_{d'}(C) \leq (d-4) + (d'-4) = g-6 < g-4$$
which contradicts relation \eqref{relationlb}.

In order to prove (i) we need to show that the case $C$ bi-elliptic also leads
to a contradiction. Let $\pi: C \lra E$ be a degree $2$ mapping onto an 
elliptic
curve $E$. Then by \cite{acgh} p.269, exercise E1, the chords $\overline{pq}$
with $p+q = \pi^* e$ for some $e \in E$ all pass through a common point $a 
\notin C$. In particular $a$ lies on a chord through all points $p_\alpha \in
\si Q_{ij}$, hence $a \in Q_{ij}$. Since $C \subset |K|^*$ is non-degenerate,
$a \in \si Q_{ij}$ and for $p_1,\ldots,p_{g-2}$ in general position
$$ a \notin \bar{D}_{ij}$$
hence $\rk Q_{ij} \leq 3$, a contradiction.

\bigskip
\noindent
It remains to show that (ii) holds. We observe that 
\begin{eqnarray*}
\rk Q_{ij} = 4 & \iff & \eta_{ij}, \nu_{ij} \in \langle \omega_i, \omega_j
\rangle \\
  & \iff & \si Q_{ij} = \bar{D}_{ij}
\end{eqnarray*}
and we conclude as before.
\end{proof}

\begin{rem}
If $C$ is trigonal or a smooth plane quintic, then for all $(g-2)$-tuples
$p_1, \ldots, p_{g-2}$ (in general position) we have
$$ \si Q_{ij} \cap C = D_{ij} \qquad \text{and} \qquad  \rk Q_{ij} = 4 $$
We can give a more precise description of the quadrics $Q_{ij}$ in both
cases: it will be enough to exhibit a set of quadrics satisfying the 
characterizing properties \eqref{charquad}.

\bigskip
\noindent
1. $C$ is {\em trigonal}. Let $|g^1_3|$ be the trigonal pencil and 
consider the complete linear series of degree $g-1$
$$ \xi = g^1_3 + D_{ij} $$
where $D_{ij}$ is as in the proof of Prop. 5.3. Then define $Q_\xi$ to be
the cone with vertex $\bar{D}_{ij} = \mathbb{P}^{g-5}$ over the smooth
quadric $|\xi|^* \times |K\xi^{-1}|^* = \mathbb{P}^1 \times \mathbb{P}^1 
\subset \mathbb{P}^3 = |K(-D_{ij})|^*$
\begin{equation} \label{defquadtri}
\begin{array}{ccc}
Q_\xi & \subset & |K|^* \\
 & & \\
\downarrow & & \downarrow^{pr} \\
 & & \\
\mathbb{P}^1 \times \mathbb{P}^1 & \map{m} & \mathbb{P}^3 
\end{array}
\end{equation}  
where $pr$ is the linear projection with center $\bar{D}_{ij}$ and $m$
is the Segre map. From this description it
is clear that the rank $4$ quadric $Q_\xi$ satisfies \eqref{charquad}.

\bigskip
\noindent
2. $C$ is a {\em smooth plane quintic} ($g=6$). Let $|g^2_5|$ be the 
associated degree $5$ linear series. We can write $D_{ij} = p_k + p_l$ for
some indices $k,l$ and, as in the trigonal case, we consider the quadric
$Q_\xi$ defined by the diagram \eqref{defquadtri} with
$$ \xi = g^2_5(-p_k) \qquad \text{and} \qquad K\xi^{-1} = g^2_5(-p_l) $$
Again we easily check that $Q_\xi$ satisfies \eqref{charquad}.
\end{rem}

\subsection{Rank $6$ quadrics and rank $2$ vector bundles}

In this section we recall a construction \cite{brivio} relating rank $6$ 
quadrics 
and rank $2$ vector bundles. We consider a rank $2$ bundle $E$ and a subspace
$V \subset H^0(E)$ which satisfy the conditions:
\begin{equation} \label{condition}
\begin{split}
 & \det E = K \\
 & \dim V = 4 \\
 & V \: \text{generates} \: E.
\end{split}
\end{equation}
\noindent
We can associate to such a bundle the following commutative diagram 
\begin{equation} \label{grasscon}
\begin{CD}
C @>\gamma>> Gr(2,V^*) \\
@VV\varphi_KV    @VVV \\
|K|^* @>\lambda^*>> \mathbb{P}(\wedge^2 V^*) = \mathbb{P}^5
\end{CD}
\end{equation}
\noindent
where $\gamma$ is the morphism (defined since we have a surjection
$\mathcal{O}_C \otimes V \ra E$)
\begin{equation} \label{defgamma}
\gamma : p \lms E_p^* \subset V^*
\end{equation}
\noindent
and $\lambda$ the map defined by taking the exterior product of global
sections of $E$
\begin{equation} \label{extprod}
\lambda : \wedge^2 V \lra H^0(\wedge^2 E) = H^0(K).
\end{equation}
\noindent
The Grassmannian $Gr(2,V^*)$ embedded in $\mathbb{P}^5$ by the Pl\"ucker 
embedding is a smooth quadric. We define $Q_{(E,V)}$ to be the inverse image of
this quadric:
\begin{equation} \label{defquad}
Q_{(E,V)} = (\lambda^*)^{-1} (Gr(2,V^*)).
\end{equation} 
Then $\rk Q_{(E,V)} \leq 6$ and  $Q_{(E,V)} \in |I(2)|$. If $h^0(E) = 4$, then
$V= H^0(E)$ and we denote the quadric $Q_{(E,V)}$ simply by $Q_E$. We have the
following lemmas:

\begin{lem}
For any pair $(E,V)$ satisfying conditions \eqref{condition}, let
$\tilde{Q}_{(E,V)}$ be the polar form of the quadric $Q_{(E,V)}$. Then
$$ \forall p,q \in C, \: p \not=q\: : \: \: \tilde{Q}_{(E,V)}(p,q) = 0 \iff
V \cap H^0(E(-p-q)) \not= \{0\} .$$
\end{lem}

\begin{proof}
We consider the dual Grassmannian $Gr(2,V)$. Now $\tilde{Q}_{(E,V)}(p,q) =0$ 
means
that the line joining the two points corresponding to the two subspaces
$H^0(E(-p)) \cap V$ and $H^0(E(-q)) \cap V$ is contained in $Gr(2,V)$. But
this is equivalent to $H^0(E(-p)) \cap H^0(E(-q)) \cap V \not= \{0\}$ and
we are done.
\end{proof}

\begin{lem}
For any pair $(E,V)$ satisfying conditions \eqref{condition}, we have
\begin{equation*}
\begin{split}
\rk Q_{(E,V)} \leq 4 \iff &  E \:  \text{contains a line subbundle} \:  L \:
\text{with} \\
& \dim H^0(L) \cap V = 2. 
\end{split}
\end{equation*}
\end{lem}

\begin{proof}
This is essentially Prop. 1.11 of \cite{brivio}. Note that if $E$ is generated
by $V$, then $\dim H^0(L) \cap V \leq 2$.
\end{proof}

\begin{rem}
We see that the definition \eqref{defquad} of the quadric $Q_{(E,V)}$ makes 
sense
even if the bundle $E$ is not generated by global sections in $V$ at a finite
number of points. We easily see that $E$ is not generated by $V$ at the point 
$p$
if and only if $p \in \mathbb{P} \ker \lambda^* \subset \si Q_{(E,V)}$. 
Moreover,
if $\rk Q_{(E,V)} \geq 5$, then we have an equality $\mathbb{P} \ker \lambda^*
= \si Q_{(E,V)}$ (see \cite{brivio} (1.9)).
\end{rem}

\bigskip

The above described construction which associates to the pair $(E,V)$ the 
quadric
$Q_{(E,V)} \in I(2)$ admits an {\em inverse} construction, i.e. we can recover 
a
bundle $E$ from a general rank $6$ quadric: consider a quadric $Q \in I(2)$ 
satisfying
\begin{equation} \label{condinvcon}
\begin{split}
 &  r = \rk Q = 5 \: \text{or} \: 6 \\
 &  \si Q \cap C = \emptyset 
\end{split}
\end{equation}
We project away from $\si Q$
$$ \delta : Q \lra \mathbb{P}^{r-1}. $$
If $r=6$, the image $\delta(Q)$ is a smooth quadric in $\mathbb{P}^5$ and 
can be realized as a Grassmannian $Gr(2,4)$. If $r=5$, $\delta(Q)$ is a linear
section of $Gr(2,4) \subset \mathbb{P}^5$. We consider the exact sequence
over $Gr(2,4)$
$$ 0 \lra \mathcal{U} \lra \mathcal{O}_{Gr}^4 \lra \bar{\mathcal{U}} \lra 0 $$
where $\mathcal{U}$ (resp. $\bar{\mathcal{U}}$) is the universal subbundle 
(resp.
quotient bundle). Since $\si Q \cap C = \emptyset$, the restriction of $\delta$
to the curve $C$ is everywhere defined and we can consider the two pairs 
(let $g = \delta_{|C}$)
\begin{equation} \label{pairs}
 (g^*\mathcal{U}^* , g^*H^0(\mathcal{U}^*)) \qquad (g^*\bar{\mathcal{U}}, 
g^*H^0(\bar{\mathcal{U}})) 
\end{equation}
which satisfy conditions \eqref{condition}. The following proposition is 
proved in
\cite{brivio} (Prop. (1.18) and (1.19))
\begin{prop}
The pairs \eqref{pairs} are the only pairs defining the quadric $Q$. If 
$\rk Q = 5$ then they are isomorphic.
\end{prop}

\begin{lem}
Consider a bundle $E$ with $h^0(E) = 4$ and satisfying \eqref{condition}.
If $\rk Q_E = 5 \: \text{or} \: 6$, then $E$ is stable
\end{lem}
\begin{proof}
Suppose that there exists a destabilizing subbundle $L \subset E$, with $\deg
L  \geq g-1$. By Lemma 5.6 we have $h^0(L) \leq 1$ and by Riemann-Roch,
$h^0(KL^{-1}) \leq h^0(L)$. But then $h^0(E) \leq h^0(L) + h^0(KL^{-1})
\leq 2$, a contradiction.
\end{proof}

\subsection{Proof of theorem 1.2}

\subsubsection{The map in \eqref{exactsequence}}

First, we prove the inclusion $\ggf \subset \ggd$. Consider a point
$p-q \in C-C$ and a point $r \in C$. The curve $t_x(C) \subset JC$, with
$x=\mathcal{O}(r + p -q)$, is contained in the fourfold $\cmc$. Therefore
a hyperplane $H$ in $|2 \Theta |$ containing $\cmc$ also contains all
tangent lines $T_{p-q}(C_x)$. Since these tangent lines
(fix $p,q$ and vary $r$) generate linearly the projectivized tangent
space $\mathbb{P}T_{p-q}JC$, we get $\pp(\mathbb{T}) \subset H$. We recall that
$\mathbb{T}$ is the linear span of the $\mathbb{P} T_{p-q}JC$ when $p,q$
vary (see section 4). Hence we get the inclusion.

\bigskip
\noindent
Now we consider the difference map
$$\gamma : C^4 \map{pr} C_2  \times C_2 \map{\phi_2} JC$$
\noindent
where $C^4$ is the 4-fold product of the curve and the first arrow $pr$
is the quotient by the transpositions $(1,2)$ and $(3,4)$ acting
on $C^4$. We denote by $\Delta_{i,j}$ the divisor in $C^4$ consisting
of 4-tuples having equal i-th and j-th entry. A straightforward
computation shows that
\begin{equation} \label{pullback2theta}
\gamma^* \mathcal{O}(2 \Theta_0) = \bigotimes_{i=1}^{4} \pi_i^* K
(- 2\Delta_{1,2} - 2\Delta_{3,4} + 2\Delta_{1,3} + 2\Delta_{1,4}
+ 2\Delta_{2,3} + 2 \Delta_{2,4})
\end{equation}
\noindent
Note that the divisor $\Delta_{1,3} + \Delta_{1,4}
+\Delta_{2,3} + \Delta_{2,4} \subset C^4$ is invariant under the
transpositions $(1,2)$ and $(3,4)$, hence comes from
an irreducible divisor in $C_2 \times C_2$, which we call
$\Delta$. We also observe that the line bundle $\pi_1^*K \otimes
\pi_2^*K(-2\Delta_{1,2})$ over $C^2$ is invariant under the
natural involution, hence comes from a line bundle $\mathcal{M}$
over $C_2$ and we have a canonical isomorphism (see e.g. \cite{brivio})
\begin{equation} \label{isoquad}
H^0(C_2,\mathcal{M}) \cong I(2).
\end{equation}
\noindent
With this notation we rewrite \eqref{pullback2theta} as
\begin{equation} \label{pullbackbis}
\phi_2^*(\mathcal{O}(2\Theta_0)) = \pi_1^* \mathcal{M} \otimes \pi_2^*
\mathcal{M}(2\Delta).
\end{equation}

\noindent
Now we want to compute the pull-back of $2\theta$-divisors
vanishing doubly on $C-C$. Let $\mathcal{J}$ be the sheaf of
ideals defining the surface $C-C \subset JC$.

\begin{lem}
If $C$ is non-trigonal, the inverse image ideal sheaf 
$\phi_2^{-1}{\cal J}\cdot \mathcal{O}_{C_2 \times C_2}$ 
is the invertible sheaf
$\mathcal{O}_{C_2 \times C_2}(-\Delta)$, hence
 $\phi_2^*\mathcal{J} =\mathcal{O}_{C_2 \times C_2}(-\Delta)$
\end{lem}

\begin{proof}
This follows from the observation that the inverse image
of $C-C$ under $\phi_2$ is isomorphic to the divisor $\Delta$.
\end{proof}

\begin{rem}
If $C$ is trigonal, the inverse image $\phi_2^{-1}(C-C)$
contains, apart from the divisor $\Delta$, a surface which is the
image of the morphism
$$
\begin{array}{rcl}
C \times C & \longrightarrow & C_2 \times C_2 \\
(p,q) & \longmapsto & (g_3^1(-p),g_3^1(-q))
\end{array}
$$
where $g_3^1$ is the trigonal series (unique if $g \geq 5$) and $g_3^1(-p)$ 
denotes the
residual pair of points in the fibre containing $p$. We deduce that 
$\phi_2^{-1}{\cal J}\cdot \mathcal{O}_{C_2 \times C_2} \subset
 \mathcal{O}_{C_2 \times C_2}(-\Delta)$ and that there exists a natural map
of $\mathcal{O}_{C_2 \times C_2}$-modules $\phi_2^*{\cal J} \lra 
 \mathcal{O}_{C_2 \times C_2}(-\Delta)$. 
\end{rem}
\bigskip
\noindent
Combining \eqref{pullbackbis} and lemma 5.10, we obtain a linear map
induced by $\phi_2$
\begin{equation} \label{defphi2}
\phi_2^* :H^0(JC, \mathcal{O}(2\Theta_0) \otimes \mathcal{J}^2) = \ggd 
\longrightarrow H^0(C_2 \times C_2,
\pi_1^* {\cal M} \otimes \pi_2^* {\cal M})
\end{equation}
This map is equivariant under the natural involutions of $JC$  and $C_2 \times
C_2$. Since all second-order theta functions are even, the
image of $\phi_2^*$ is contained in ${\rm Sym}^2
H^0(C_2, {\cal M}) \cong \sym^2 I(2)$, by \eqref{isoquad}.
To summarize, we have shown that
$$ \ggf = \ker( \phi_2^* : \ggd \lra \sym^2 I(2)). $$

\subsubsection{Surjectivity of $\phi_2^*$}
The key point of the proof is the following proposition
\begin{prop}
We have a commutative diagram
\begin{equation} \label{comdiagram}
\begin{CD}
\mathcal{W}(4) @>D>> \mathbb{P} \ggd \\
@VVQV @VV\phi_2^*V \\
|I(2)| @>Ver>> \mathbb{P} \sym^2 I(2)              
\end{CD}
\end{equation}
where the notation is as follows:

$\mathcal{W}(4) = \{ [E] \in \modk \: | \: \dim H^0(E) = 4 \:
\text{and} \: E \: \text{globally generated} \}$

$Q$ is the map described in section 5.2

$Ver$ is the Veronese map 

$D$ is the map \eqref{defd}
\end{prop}

\begin{proof}
Consider a bundle $E \in \mathcal{W}(4)$. By \cite{laszlo} Prop. V.2., we
have an inequality
$$ \mult_{p-q}(D(E)) \geq h^0(E(p-q)) \geq 2  $$
hence we see that $D(E) \in \mathbb{P} \ggd$. To show commutativity, it is
enough to check that the zero divisors of the two sections
$$ Q_E \otimes Q_E \qquad \phi_2^*(D(E)) $$
(regarded as sections of $\pi_1^* \mathcal{M} \otimes  \pi_2^* \mathcal{M}$
over $C_2 \times C_2$ \eqref{defphi2}) coincide as sets. Hence, by Lemma
5.5 and \eqref{defd}, it is enough to show the following equivalence: for any
four {\em distinct} points $p,q,r,s \in C$
$$ h^0(E(-p-q)) > 0 \: \: \text{or} \: \: h^0(E(-r-s)) >0 \: \: \iff \: \:
h^0(E(p+q-r-s)) > 0 $$
The $\Rightarrow$ implication is obvious ($D(E)$ is symmetric). To prove the
$\Leftarrow$ implication, we suppose that $h^0(E(-p-q)) = h^0(E(-r-s)) = 0$.
Then, by Riemann-Roch and Serre duality, we have $h^0(E(p+q))=4$. Since
$h^0(E) = 4$, we see that all global sections of $E(p+q)$ vanish at the
points $p$ and $q$. Supposing that there exists a non-zero section $\varphi$
of $E(p+q-r-s)$, then $\varphi$ vanishes at $p,q$, contradicting
$h^0(E(-r-s))=0$, hence $h^0(E(p+q-r-s)) = 0$.  
\end{proof}

 From now on, we assume that $C$ is {\em non-trigonal}.
We consider $g$ points $p_1,\ldots,p_g$ in general position and their
associated Petri quadrics $Q_{ij}$ for $1 \leq i,j \leq g-2$ (section 5.1),
which form a basis of $I(2)$. Then in order to prove surjectivity of
$\phi_2^*$, it is enough, by Prop. 5.12, to construct a set of $\binom{h+1}{2}$
vector bundles (with $h = \binom{g-2}{2} = \dim I(2)$) in $\mathcal{W}(4)$,
which generate linearly $\sym^2 I(2)$. First we suppose that $C$ is not a 
smooth plane quintic. We proceed in 3 steps.

\bigskip
\noindent
{\em Step 1}

\bigskip
\noindent
By Prop. 5.3, the
Petri quadrics $Q_{ij}$ satisfy conditions \eqref{condinvcon}, so
(Prop. 5.8) we can construct for each $i,j$ two pairs (see \eqref{pairs}) of
bundles $(E_{ij}^{(1)}, V^{(1)})$ and  $(E_{ij}^{(2)}, V^{(2)})$, which
define the quadric $Q_{ij}$.

\begin{lem}
For general points $p_1, \ldots, p_{g-2}$, the bundles $E_{ij}^{(1)}, 
E_{ij}^{(2)}$
are stable, distinct and $h^0(E_{ij}^{(1)}) = h^0(E_{ij}^{(2)}) = 4$. 
\end{lem} 

\begin{proof}
We can give a different description  of the bundles $E_{ij}^{(1)}$ and 
$E_{ij}^{(2)}$ using extension spaces. Let $D = D_{ij} + p_i$ and consider
extensions of the form
$$ 0 \lra \mathcal{O}(D) \lra \mathcal{E}_\varepsilon \lra 
\mathcal{O}(K-D)  \lra 0 \hspace{1cm} (\varepsilon) $$
These extensions are classified by an extension class $\varepsilon \in
|2K-2D|^* = \mathbb{P}^{g+2}$. Since $h^0(D) = 1$ and $h^0(K-D) = 3$, we
see that $h^0(\mathcal{E}_\varepsilon) = 4$ if and only if $\varepsilon
\in \ker m^* = (\coker m)^*$, where $m$ is the multiplication map
$$ m: \sym^2 H^0(K-D) \lra H^0(2K-2D) $$
which is injective for general points $p_i$. Note that $\dim \coker m =
g-3$. Furthermore, consider a point $p_\alpha \in D_{ij}$ and the
multiplication map (which is injective)
$$m_\alpha : H^0(K-D+p_j+p_\alpha) \otimes H^0(K-D-p_j-p_\alpha) \lra
H^0(2K-2D).$$
Then $h^0(\mathcal{E}_\varepsilon(-p_j-p_\alpha)) > 0 \iff \varepsilon
\in \ker m^*$. We observe that the image $\mathrm{im} \: m_\alpha \subset
H^0(2K-2D)$ under the canonical surjection $H^0(2K-2D) \lra \coker m$ is
a one-dimensional subspace, which we denote by $Z_\alpha$. Consider 
now a hyperplane $H$ in $\coker m$, which contains the linear span of the
$Z_\alpha$ for $\alpha$ such that $p_\alpha \in D_{ij}$ (we will see
a posteriori that such an $H$ is unique, for dimensional reasons), so that we
obtain an extension class $\varepsilon = \varepsilon(H) \in 
\mathbb{P} (\coker m)^* \subset |2K-2D|^*$. By construction, we have 
$h^0(\mathcal{E}_\varepsilon) = 4$ and $h^0(\mathcal{E}_\varepsilon(-p_\alpha
-p_\beta)) >0$ if $\{ \alpha, \beta \} \not= \{ i,j \}$, hence, by 
\eqref{charquad} and lemma 5.5, we get $Q_{\mathcal{E}_\varepsilon} =
Q_{ij}$, and $\mathcal{E}_\varepsilon = E_{ij}^{(1)}$.

The other bundle $E_{ij}^{(2)}$ defining the quadric $Q_{ij}$ is constructed
in the same way using the divisor $D' = D_{ij} +p_j$ (instead of $D$).
Then we have $E_{ij}^{(1)} \not= E_{ij}^{(2)}$. Indeed, an isomorphism
$E_{ij}^{(1)} \rightarrow E_{ij}^{(2)}$ would imply the existence
of a nonzero section of $\mathcal{H}om(\mathcal{O}(D), \mathcal{O}(K-D')) =
\mathcal{O}(K-D-D')$, but then the points $p_1,\ldots,p_{g-2}$ are not 
in general position.

Finally, stability follows from Lemma 5.9
\end{proof}
We deduce from this lemma and Prop. 5.12 that $Q_{ij} \otimes Q_{ij} \in 
\mathrm{im}\: 
\phi_2^*$.

\bigskip
\noindent
{\em Step 2}

\bigskip
\noindent
Consider three distinct indices $i,j,k$. Then all quadrics of the pencil
\begin{equation} \label{pencil}
(\lambda Q_{ij} + \mu Q_{ik})_{\lambda, \mu \in \mathbb{C}}
\end{equation}
have rank less than or equal to $6$. This follows from expression 
\eqref{Petriquad}
of Petri's quadrics, namely
$$ \lambda Q_{ij} + \mu Q_{ik} = \omega_i (\lambda \omega_j + \mu \omega_k) -
\bar{\eta} \omega_{g-1} - \bar{\nu} \omega_g - \bar{b} \omega_{g-1}
\omega_g $$
with $\bar{\eta} = \lambda \eta_{ij} + \mu \eta_{ik}, \bar{\nu} = \lambda
\nu_{ij} + \mu \nu_{ik}, \bar{b} = \lambda b_{ij} + \mu b_{ik}$. Now a
general element of the pencil \eqref{pencil} satisfies conditions
\eqref{condinvcon}, since these are open conditions and are satisfied by
$Q_{ij}$ and $Q_{ik}$. Again by openness and Lemma 5.13, it follows
that the two bundles associated with such a general quadric have $4$
sections and are stable. Let us pick such a bundle $E_{ijk}$ defining
the quadric $\lambda_0 Q_{ij} + \mu_0 Q_{ik}$, for $\lambda_0, \mu_0
\not= 0$. Then we have in $\sym^2 I(2)$
$$\phi_2^* D(E_{ijk}) = \lambda_0^2 Q_{ij} \otimes Q_{ij} + \mu_0^2
Q_{ik} \otimes Q_{ik} + 2\lambda_0 \mu_0 Q_{ij} \otimes Q_{ik} $$
hence, $Q_{ij} \otimes Q_{ik} \in \mathrm{im}\: \phi_2^*$.

\bigskip
\noindent
{\em Step 3}

\bigskip
\noindent
Consider four distinct indices $i,j,k,l$. Then all quadrics of the
$2$-dimensional family $\mathcal{F}_{(ij)(kl)} = \mathbb{P}^1 \times
\mathbb{P}^1$ (here $(\lambda,\lambda'), (\mu,\mu')$ are a set of
homogeneous coordinates) are given by an expression:
\begin{eqnarray*}
\mu \lambda Q_{ik} + \mu \lambda' Q_{il} + \mu' \lambda Q_{ik} +
\mu' \lambda' Q_{jl} = \\
 (\mu \omega_i + \mu' \omega_j)(\lambda \omega_k + \lambda' \omega_l)
 - \bar{\eta} \omega_{g-1} - \bar{\nu} \omega_g - \bar{b} \omega_{g-1}
\omega_g,
\end{eqnarray*}
where $\bar{\eta},\bar{\nu},\bar{b}$ depend on 
$(\lambda,\lambda'), (\mu,\mu')$,
see \eqref{Petriquad}. The same holds for the two
families obtained by permuting indices
$$ \mathcal{F}_{(ik)(il)} \: : \: 
\mu \lambda Q_{ij} + \mu \lambda' Q_{il} + \mu' \lambda Q_{jk} +
\mu' \lambda' Q_{kl} $$
$$ \mathcal{F}_{(il)(kj)} \: : \: 
\mu \lambda Q_{ik} + \mu \lambda' Q_{ij} + \mu' \lambda Q_{lk} +
\mu' \lambda' Q_{jl} $$
As in step 2, we see that a general member of these 3 families satisfies
\eqref{condinvcon} and we can pick 3 stable vector bundles 
$E_{(ij)(kl)}, E_{(ik)(jl)}, E_{(il)(kj)}$ with 4 sections defining
the 3 quadrics in these families with coordinates $(\lambda_0,\lambda'_0)
(\mu_0,\mu'_0)$ for some $\lambda_0, \lambda'_0, \mu_0, \mu'_0 \not= 0$.
Now we can write in $\sym^2 I(2)$
$$ \phi_2^* D(E_{(ij)(kl)}) = 2 \mu_0 \mu'_0 \lambda_0 \lambda'_0 
(Q_{ik} \otimes Q_{jl} + Q_{il} \otimes Q_{jk}) + \alpha $$
$$ \phi_2^* D(E_{(ik)(jl)}) = 2 \mu_0 \mu'_0 \lambda_0 \lambda'_0 
(Q_{ij} \otimes Q_{kl} + Q_{il} \otimes Q_{jk})  + \beta $$
$$ \phi_2^* D(E_{(il)(kj)}) = 2 \mu_0 \mu'_0 \lambda_0 \lambda'_0 
(Q_{ik} \otimes Q_{jl} + Q_{ij} \otimes Q_{kl}) + \gamma $$
for some $\alpha, \beta, \gamma \in \mathrm{im} \: \phi_2^*$ (see step 1 and
2). But these linear equations immediately imply that the three
symmetric tensors $Q_{ij} \otimes Q_{kl}, Q_{ik} \otimes Q_{jl},
Q_{il} \otimes Q_{jk} \in \mathrm{im} \: \phi_2^*$ and we are done.

\bigskip

To complete the proof we need to consider the case when $C$ is a smooth
plane quintic ($g=6$). We will show that the map $Q$ in diagram
\eqref{comdiagram} is dominant. By \cite{arba} Prop. 3.2, the locus of 
rank 4 quadrics is a cubic hypersurface in $|I(2)| = \mathbb{P}^5$ and a
general quadric has rank $6$ (i.e. is smooth). Consider any smooth quadric 
$Q \in |I(2)|$ and one of the associated pairs $(E,V)$ defining $Q$ 
\eqref{pairs}. It will be sufficient to show that $h^0(E)= 4$, hence
$E \in \mathcal{W}(4)$. By \cite{prev} Thm. 8.1 (3), $h^0(E) \geq 5$ if and
only if $E$ is an extension of the form
$$ 0 \lra g^2_5 \lra E \lra g^2_5 \lra 0. $$ 
From this we see that if $h^0(E) = 5$, $\dim V \cap H^0(g_5^2) \geq 2$,
hence by Lemma 5.6 $\rk Q \leq 4$ and if $h^0(E)= 6$, then $E = g^2_5 \oplus
g^2_5$ and we can also conclude that $\rk Q \leq 4$, contradicting $Q$
smooth.

\subsection{Another proof of a theorem by M. Green}

As a consequence of
Thms. 1.1 and 1.2 we get another proof of the 
following theorem (in the case of non-trigonal curves) due to M. Green
(\cite{green}, see also \cite{smith})
\begin{thm}
For $C$ non-trigonal, the projectivized tangent cones at singular points 
of $\Theta$ span $I(2)$.
\end{thm}

\begin{proof}
For all $\xi \in \si \Theta$ with $h^0(\xi) = 2$, the split bundle
$E = \xi \oplus K\xi^{-1} \in \mathcal{W}(4)$. Then the associated quadric 
$Q_E$
has rank $4$ and can be described as a cone over the smooth quadric
$\mathbb{P}^1 \times \mathbb{P}^1 \subset \mathbb{P} V^*  = \mathbb{P}^3$
(as in diagram \eqref{defquadtri}) where $V$ is the image of the
multiplication map 
$$H^0(\xi) \otimes H^0(K\xi^{-1}) \lra H^0(K). $$
Then $Q_E$ is the 
projectivized tangent cone at $\xi \in \si \Theta$.

Suppose now that the image under $Q$ of $\si \Theta$ in $|I(2)|$ is 
degenerate. By Prop. 5.12 and Thm. 1.1 (2), we see that the image of 
$\phi_2^*$ is also degenerate, contradicting Thm. 1.2.
\end{proof}

\section{The space $\ggg$}

\subsection{Proof of theorem 1.3}

In this section we regard $2\theta$-functions as holomorphic functions
on $\cc^g$. For the proof of Theorem 1.3 we  need the following lemma:

\begin{lem}
The following statements are equivalent
\begin{enumerate}
\item[(1)]
$f \in \ggg$
\item[(2)]
For all $a_1,a_2,a_3,a_4 \in \ca$
\begin{equation*}
\begin{split} 
(q_{12}q_{34})^4 f(w(a_1 + a_2 - a_3 -a_4)) + (q_{14}q_{23})^4
 f(w(a_1 &+ a_4 - a_2 -a_3)) + \\ &(q_{13}q_{24})^4
f(w(a_1 +a_3 -a_2 -a_4)) = 0
\end{split}
\end{equation*}
\end{enumerate}
\end{lem}

\begin{proof}
We will derive identity (2) from Prop. 3.3. First we observe that
the left-hand side $PD_{a_1}D_{a_2}D_{a_3}D_{a_4}\theta_2(0)$ of Prop. 3.3 (1)
is symmetric in the four variables. In particular it is symmetric in $a_2, 
a_4$,
which leads to the equality
\begin{equation} \label{eee}
\tau_{1234} + \tau_{2314} = \tau_{1432} + \tau_{4312} + 2 \sigma_{1234}
+ 2 \sigma_{2134} + 2 \sigma_{3124} - 2 \sigma_{1432} - 2 \sigma_{4132}
 - 2 \sigma_{3142}
\end{equation}
Combining \eqref{eee} and prop.3.3 (1) we can write
$$ \frac{1}{2} PD_{a_1}D_{a_2}D_{a_3}D_{a_4}\theta_2(0) = \tau_{1324}
+ \tau_{1432} + \tau_{4312} + \text{some\: } \sigma 's$$
Using Prop. 3.3 (2) and the two relations obtained from it
by interchanging $a_1$ with $a_3$ and $a_1$ with $a_4$, we can write

$$\frac{1}{2} PD_{a_1}D_{a_2}D_{a_3}D_{a_4}\theta_2(0) =  
 {\left( \frac{q_{12}q_{34}} 
{q_{13}q_{14}q_{23}q_{24}} \right)}^2 P\theta_2(w(a_1+a_2-a_3-a_4)) +$$
$$  {\left( \frac{q_{32}q_{14}} 
{q_{13}q_{34}q_{21}q_{24}} \right)}^2 P\theta_2(w(a_1+a_4-a_3-a_2)) + 
 {\left( \frac{q_{24}q_{31}} 
{q_{43}q_{14}q_{23}q_{21}} \right)}^2 P\theta_2(w(a_1+a_3-a_2-a_4)) - X$$
where $X$ is the following sum of $\sigma$'s
$$ ( \sigma_{1243} + \sigma_{1234} + \sigma_{1432} ) +
 ( \sigma_{3142} + \sigma_{3241} + \sigma_{3214} ) +
 ( \sigma_{4132} + \sigma_{4213} + \sigma_{4231} ) $$
But the three terms within each pair of parentheses add up to zero (use the 
definition of
$\sigma$ \eqref{defsigma} and the fact that $\xi$ is skew-symmetric in its 
variables), hence
$X=0$ and we are done. 
\end{proof}

We are now in a position to prove Theorem 1.3. First we show the inclusion
$\ggg \subset \ggd$. We fix $f \in \ggg$ and three points $a_2,a_3,a_4 \in
\ca$. We consider $a_1$ as a canonical coordinate and derive two (resp. three) 
times
with respect to $a_1$ and take the value at the point $a_1 =a_4$. This way,
we obtain two equations among vectors in $\cc^{2^g}/\tankum$
$$ D_{a_4}^2 \theta_2(w(a_2-a_3)) + D_{a_4} \theta_2(w(a_2-a_3)) 
\left[4\partial \log \frac{q_{42}}{q_{43}} \right] = 0 $$
$$ D_{a_4}^2 \theta_2(w(a_2-a_3)) \left[ \partial \log q_{42} q_{43} \right]
+ D_{a_4} \theta_2(w(a_2-a_3)) \left[ \partial^2 \log \frac{q_{42}}{q_{43}}
  + 4 \partial \log q_{42}q_{43} \cdot \partial \log \frac{q_{42}}{q_{43}}
\right] = 0 $$
where $\partial$ means derivative with respect to the first variable of the
prime form $q$. Hence we get, for all $a_2,a_3,a_4 \in \ca$, a system of two
linear equations involving the vectors $D_{a_4}^2\theta_2(w(a_2 - a_3))$ and
$D_{a_4}\theta_2(w(a_2-a_3))$ whose determinant
$$\partial^2 \log \frac{q_{42}}{q_{43}}$$
is non-zero on an open subset of $\ca^3$. Hence the two vectors
$D_{a_4}^2\theta_2(w(a_2 - a_3))$ and $D_{a_4}\theta_2(w(a_2-a_3))$ are
zero on an open subset of $\ca^3$, so they are identically zero.
This implies that $f \in \ggd$.

\bigskip

The inclusion $\ggf \subset \ggg$ and the second 
statement of Theorem 1.3 can easily be deduced from the 
commutativity of the diagram
\begin{equation} \label{codi}
\begin{array}{ccc}
\spansing = \ggd & \map{\phi_2^*} & \sym^2 I(2) \\
 & & \\
\downarrow^\alpha & \swarrow & \\
 & & \\
I(4) & & 
\end{array}
\end{equation}
where $\alpha$ is the map which associates to a $2\theta$-divisor its  
projectivized quartic tangent cone at the origin (i.e. the degree 4 term of 
the Taylor expansion at the origin), $\phi_2^*$ is as in \eqref{defphi2} and 
the
diagonal arrow is the multiplication map $m$. By definition we have
$$ \ker \alpha = \ggg \qquad \ker \phi_2^* = \ggf $$
To check the commutativity of \eqref{codi}, by Theorem 1.1 (2) it is 
enough to
check that $\alpha(D(E)) = m(\phi_2^*(D(E)))$ for the bundle $E = \xi \oplus
K\xi^{-1}$, with $\xi \in \si \Theta$ and $h^0(\xi) = 2$, i.e. that $D(E) = 
\Theta_\xi + 
\Theta_{K\xi^{-1}}$. This follows from Prop. 5.12 and the description
of $Q_E$ (see proof of Theorem 5.14).

\begin{rem}
We don't know how to find a general formula for
$\dim \ggg$. In the examples that follow, we give 
$\dim \ggg$ for $g \leq 7$.
\end{rem}  

\subsection{Examples}

\subsubsection{Curves of genus less than $6$}

For any non-hyperelliptic curve of genus $g \leq 5$, we have 
$$ \dim \ggg = 0. $$
This is an easy consequence of \eqref{codimgg} and Thms. 1.1, 1.2 and 1.3 
(For a trigonal genus 
$5$ curve, we also use Prop. 7.2).

\subsubsection{Curves of genus $6$}

Consider first a genus $6$ curve which is not trigonal or a smooth
plane quintic. In order to determine the quadratic syzygies in $I(2)$, i.e.
$\ker m$, we consider the rational map induced by the linear system
$|I(2)|$ on the canonical space $|K|^* = \pp^5$
$$ \pi : |K|^* \lra |I(2)|^* = \pp^5. $$
Note that $\pi$ is defined away from the canonical curve. We need
the following lemma (due to S. Verra).

\begin{lem}
The rational map $\pi$ is finite of degree $2$.
\end{lem}

\begin{proof}
Let $H$ be a general hyperplane in $I(2)$ and choose a basis of
quadrics $\{ Q_0,Q_1,\ldots,Q_4\}$ of $H$. By a Bertini argument, 
the intersection in $|K|^*$ of three general quadrics $Q_i \in |I(2)|$
(for $i=0,1,2$) is a smooth, irreducible surface $S$ (hence a $K3$-surface),
which contains the canonical curve. Then the linear equivalence class of
the two divisors determined by the quadrics $Q_3,Q_4$ on the surface
$S$ is $2h-c$, where $h$ is the hyperplane section and $c$ the class
of the canonical curve. We compute
$$ (2h-c)^2 = 4h^2 - 4h\cdot c + c^2 = 2 $$
since $h^2 = 8, h\cdot c = 2g-2 = 10$ and $c^2 = 10$ (by the adjunction
formula with $\omega_S = \mathcal{O}_S$). Therefore the fibre of $\pi$ over the
point determined by $H$ consists of $2$ points and we are done. 
\end{proof}

\bigskip

It follows from Lemma 6.3 that $\pi$ is onto, hence $\ggg = \ggf$. By 
Thm. 8.1 (1) \cite{prev} and Thm. 5.1 (1) \cite{mukai2}, there exists a unique
stable bundle $E_{max} \in \modk$ with maximal number of sections
$h^0(E_{max}) = 5$. Since $C_2 - C_2 \subset D(E_{max})$, we have 
$D(E_{max}) \in \pp \ggf$. Moreover by \eqref{codimgg} and Thms. 1.1 and 1.2, 
we have 
$\dim \ggf = 1$, from which we deduce that 
$$ \pp \ggg = \pp \ggf = D(E_{max}). $$

As a complement to the examples of Brill-Noether loci of $\modk$ provided
in \cite{prev} we add a geometric description of the divisor
$D(E_{max}) \subset JC$. Let $L = g^1_4$ be a tetragonal series on $C$
and $\varphi_L$ be the associated surjective Abel-Jacobi map
\begin{eqnarray*}
\varphi_L : C_6 & \lra & JC \\
  (p_1,\ldots,p_6) & \lms & K^{-1}L(p_1+\ldots + p_6)
\end{eqnarray*}
Let $\pi_L$ be the map induced by the base point free linear series
$|KL^{-1}|$
$$ \pi_L : C \lra |KL^{-1}|^* = \pp^2.$$
If $C$ is bi-elliptic, the image is a smooth plane cubic. Otherwise 
(general case) $\pi_L$ maps $C$ birationally to a nodal plane sextic. 

\begin{prop}
Let $\mathcal{S} \subset C_6$ be the divisor consisting of sextuples
$(p_1,\ldots,p_6)$ such that the points $\pi_L(p_i)$ lie on a conic in
$|KL^{-1}|^* = \pp^2$. Then
$$ \varphi_L(\mathcal{S}) = D(E_{max}) $$
\end{prop}

\begin{proof}
By \cite{mukai2} Prop. 3.1 (see also \cite{prev} example 3.4), we 
can write $E_{max}$ as an extension
\begin{equation} \label{esemax}
0 \lra L \lra E_{max} \lra KL^{-1} \lra 0 
\end{equation}
By definition, we have $\lambda \in D(E_{max}) \iff h^0(E_{max} \otimes
\lambda) > 0$. First, we see that if $h^0(L \lambda) > 0$ or
$h^0(L \lambda^{-1}) > 0$, then the exact sequence \eqref{esemax}
implies that $\lambda \in D(E_{max})$ (note that $D(E_{max})$ is
symmetric). Therefore we can assume that $h^0(L \lambda) =
h^0(L \lambda^{-1}) = 0$ or equivalently that
$h^0(KL^{-1}\lambda) = h^0(KL^{-1} \lambda^{-1})
=1$. Writing  the long exact sequence associated to \eqref{esemax},
we see that $H^0(E_{max} \otimes \lambda) = \ker (\delta :
H^0(KL^{-1} \lambda) \lra H^1(L\lambda))$. Therefore 
$h^0(E_{max} \otimes \lambda) >0$ if and only if the image of the
multiplication map 
$$H^0(KL^{-1}\lambda) \otimes H^0(KL^{-1}\lambda^{-1}) \lra
H^0(K^2L^{-2}) $$
is contained in the hyperplane $\sym^2 H^0(KL^{-1}) \subset H^0(K^2L^{-2})$
which defines the extension class of $E_{max}$ in $|K^2L^{-2}|^*$ (see
\cite{prev} example 3.4). Since $h^0(KL^{-1}\lambda) = 1$, there exists
a unique sextuple $(p_1,\ldots,p_6) \in C_6$ such that $\lambda = 
\varphi_L((p_i))$. Then we deduce from the commutative diagram
$$
\begin{CD}
C @>\varphi_{K^2L^{-2}}>> |K^2L^{-2}|^* = \pp^6 \\
@VV\pi_LV @VVprV \\
|KL^{-1}|^* = \pp^2 @>Ver>> \pp \sym^2 H^0(KL^{-1})^* = \pp^5
\end{CD}
$$
that $h^0(E_{max} \otimes \lambda) > 0$ if and only if the $6$ points
$\varphi_{K^2L^{-2}}(p_i)$ lie on a hyperplane in $\pp^6$ which is the
inverse image under the projection map $pr$ of a hyperplane in
$\pp^5$. But this last condition says that the $6$ points $\pi_L(p_i)$
lie on a conic in $|KL^{-1}|^*$. Finally, we notice that if 
$h^0(L\lambda^{-1}) >0$ ($\iff h^0(\mathcal{O}(\sum p_i)) \geq 2$)
then there exists a divisor $D = \sum q_i$ in the linear system
$|\sum p_i|$ such that the $\pi_L(q_i)$ lie on a conic (e.g. we can 
show that $\mathcal{S}$ is an ample divisor on $C_6$).
\end{proof}

\bigskip

Consider now the case of a {\em smooth plane quintic}. By Thm. 8.1 (3) 
\cite{prev}
there exists a unique $S$-equivalence class $\epsilon_{max}$ such that
if $h^0(E) \geq 5$, then $[E] = \epsilon_{max}$. In particular, 
$\epsilon_{max} = [g^2_5 \oplus g^2_5 ]$. Moreover, using an explicit basis
of quadrics of $I(2)$, we can show that the map $\pi$ is birational (for
more details see \cite{prev} Thm. 5.5), hence $\ker m = \{ 0\}$. As before,
we deduce from Thms. 1.1 and 1.2 that 
$$\pp \ggg = \pp \ggf = D(\epsilon_{max}) = 2\Theta_{g^2_5}. $$

\bigskip

Consider now the case of a {\em trigonal curve}. By Thm. 8.1 (2) \cite{prev}
there exists a projective line $\pp^1_{bun}$ of stable bundles $E$ with 
$h^0(E) = 5$,
hence $\pp^1_{bun} \subset \pp \ggf$. Using Prop. 7.2, we compute that
$\dim \ggf = 2$ and (using section 7) that $\dim \ker m = 1$. Hence we have
$$ \pp^1_{bun} = \pp \ggf \qquad \dim \ggg = 3 $$

\subsubsection{Curves of genus $7$}

For a non-tetragonal genus $7$ curve, we have $\dim \ggg/\ggf = 1$ (by
\cite{mukai} thm 4.2). So we obtain $\dim \ggf = 9$ and $\dim \ggg = 10$.

\section{Trigonal curves}

Let $C$ be a trigonal curve with $g \geq 5$ and $g^1_3$
its unique trigonal series. By remark 5.11, we obtain as in \eqref{defphi2} a
linear map $\phi_2^* : \ggd \lra \sym^2 I(2) $. The aim of this section is
to compute the rank of $\phi_2^*$ (cf. section 5.3.2).
\bigskip

First we need to quote some results about quadrics containing a rational
normal scroll from \cite{arba}, which we state here for the case of the
degree $g-2$ surface $S \subset |K|^*$ ruled by the pencil of trisecants 
to the canonical curve. For a trigonal curve, the space $I(2)$ and the
space of quadrics $I_S(2)$ containing the surface $S$ are isomorphic.
Let $V = H^0(C,K-g^1_3)$. We choose two sections $s_0,s_1 \in H^0(g^1_3)$
and consider the isomorphism $\beta_0$ (resp. $\beta_1$) induced by  
multiplication by the section $s_0$ (resp. $s_1$)
$$ \beta_0 : V \map{\sim} V_0 \subset H^0(K) \qquad 
 \beta_1 : V \map{\sim} V_1 \subset H^0(K) $$
where $V_i = H^0(C, K-D_i)$ and $D_i$ is the zero divisor of the section
$s_i$. We then define a linear map
$$ \beta : \Lambda^2 V \lra \sym^2 H^0(K) $$
by setting
\begin{equation}\label{defbeta}
\beta(v\wedge w) = \beta_0(v) \otimes \beta_1(w) - \beta_0(w) \otimes
\beta_1(v)
\end{equation}
which is a quadric of rank less than or equal to 4 containing $S$. One 
checks that $\beta(v\wedge w)$ does not depend on the choice of the
sections $s_0, s_1$. Then Prop. 2.14 \cite{arba} says that $\beta$ induces
an isomorphism 
\begin{equation} \label{isobeta}
\beta : \Lambda^2 V  \map{\sim} I_S(2)
\end{equation}

\bigskip
We also define a rational map $\delta : \mathcal{W}(4) \lra Gr(2,V)$ as 
follows. Consider a semi-stable bundle $E \in \mathcal{W}(4)$ (see
\eqref{comdiagram}). Then, by \cite{mukai2}, Prop. 3.1, $h^0(E(-g^1_3)) \geq
1$; it can be shown that $h^0(E(-g^1_3))=1$
for a general bundle $E \in \mathcal{W}(4)$,
i.e. $E$ can be uniquely written as an extension
\begin{equation} \label{exte}
0 \lra g^1_3 \lra E \map{\pi} K-g^1_3 \lra 0
\end{equation}
and we define $\delta (E) = {\rm im} \  (H^0(E) \map{H^0(\pi)} H^0(K-g^1_3) =
V)$. Then we can prove the following

\begin{lem}
The map $Q: \mathcal{W}(4) \lra |I(2)|$ defined in section 5.2 factorizes as
follows
$$ \mathcal{W}(4) \map{\delta} Gr(2,V) \map{Pl} \pp(\Lambda^2 V) \map{\beta}
|I_S(2)| = |I(2)| $$
where $Pl$ is the Pl\"ucker embedding of the Grassmannian.
\end{lem}
\begin{proof}
Consider $E \in \mathcal{W}(4)$ with $h^0(E(-g^1_3)) = 1$ and identify
$H^0(g^1_3)$ with a $2$-dimensional subspace of $H^0(E)$. We choose a basis
$\{ s_0,s_1 \}$ of $H^0(g^1_3)$. Let $R = \beta \circ Pl \circ \delta (E)
\in |I(2)|$ be the associated quadric. In order to show that $R = Q_E$ it
is enough to show that their associated polar forms, which we view as 
global sections of the line bundle $\mathcal{M}$ over $C_2$ (see 
\eqref{isoquad}), coincide. Hence it is enough to show the implication
$$\forall p,q \in C,\:  p\not= q,\: h^0(g^1_3(-p-q))=0 \hspace{15mm}
\tilde{Q}_E(p,q) =0 \: \Longrightarrow  \tilde{R}(p,q) = 0 $$
But by Lemma 5.5, the assumption $\tilde{Q}_E(p,q) = 0$ means that 
there exists a section $a \in H^0(E)$ vanishing at $p$ and $q$. Since
$h^0(g^1_3(-p-q)) = 0$, we have $a \notin H^0(g^1_3)$ and we can find
a section $b \in H^0(E)$ such that $\{ s_0,s_1,a,b \}$ is a basis
of $H^0(E)$. Then $H^0(\pi)$ induces a linear isomorphism $\cc a \oplus
\cc b \map{\sim} \delta (E) \subset V$. Let $u,v \in \delta (E)$ be the
images of $a,b$ under this isomorphism. Then we see that 
$$ \beta_0(u) = s_0 \wedge a \in H^0(K) \qquad  
\beta_0(v) = s_0 \wedge b \in H^0(K). $$
The same holds for $\beta_1$ and $s_1$. By \eqref{defbeta} we have
\begin{equation*}
\begin{split}
2\tilde{R}(p,q) = (s_0 \wedge a)(p)\cdot (s_1 \wedge b)(q) &+
 (s_0 \wedge a)(q)\cdot (s_1 \wedge b)(p) \\
 &- (s_0 \wedge b)(p)\cdot (s_1 \wedge a)(q) -
 (s_0 \wedge b)(q)\cdot (s_1 \wedge a)(p)
\end{split}
\end{equation*}
and this expression, which does not depend on the choice of the basis
$\{ s_0,s_1,a,b \}$, is obviously zero if $a \in H^0(E(-p-q))$.
\end{proof}

We consider now the commutative diagram \eqref{comdiagram}. 
Lemma 7.1 implies that the inverse image under the Veronese map of a
hyperplane in $\sym^2 I(2)$ containing ${\rm im} \: \phi_2^*$ is a 
quadric in $|I(2)|$ containing the Grassmannian $Gr(2,V)$. Moreover
any such quadric comes from an element in the annihilator
$({\rm im} \: \phi_2^*)^\perp$. Hence we obtain an isomorphism
\begin{equation} \label{isogrphi}
I_{Gr(2,V)}(2) \cong ({\rm im} \: \phi_2^*)^\perp .
\end{equation}
But the degree $2$ part $I_{Gr(2,V)}(2)$ of the ideal of the Grassmannian 
$Gr(2,V)$ is
isomorphic to the vector space $\Lambda^4 V$ generated by the
Pl\"ucker equations (see e.g. \cite{mukai2}). Hence 
we have shown

\begin{prop}
The corank of $\phi_2^*$ is $\binom{g-2}{4}$.
\end{prop}

\section{Concluding remarks}

\noindent
{\bf 1.} Some analytic versions of Thms. 1.1 and 1.2 were proved by 
Gunning in the case of a  
non-hyperelliptic curve of genus less than $6$ (\cite{gu2} Thm.
6). Furthermore, the statement of Thm. 1.3 was 
proposed as plausible
in \cite{gu1},
p. 70.

\bigskip
\noindent
{\bf 2.} It is natural to ask whether  
the three main theorems can be generalized to
analogous subseries. For this purpose, we introduce the following
subspaces of $\gg$
$$
\begin{array}{rcll}
\Gamma_{[n]\cdot 0} & = & \{ D \: | \: \mult_0(D) \geq 2n \} & 
\text{for} \: n \geq 2 \\
 & & & \\
\Gamma_{dd} & = & \{ D \: | \: C_{d+1} - C_{d+1} \subset D \} &
\text{for} \: d \geq 0 \\
 & & & \\
\Gamma_{dd}^{(2)} & = & \{ D \: | \: \mult_{C_{d+1} - C_{d+1}}(D) \geq 2 \} &
\text{for} \: d \geq 0 
\end{array}
$$
where $C_{d+1} - C_{d+1}$ is the image of the difference map ($d \geq 0$)
\begin{eqnarray*}
\phi_{d+1} : C_{d+1} \times C_{d+1} & \lra & JC \\
(D,D') & \lms & \mathcal{O}(D-D')
\end{eqnarray*}
The following  inclusions are obvious
\begin{equation} \label{obviincl}
\Gamma_{(d+1)(d+1)} \subset \Gamma_{dd}^{(2)} \subset \Gamma_{dd}
\end{equation}
and one might expect that the following holds (see Thm. 1.3)
\begin{equation} \label{serincl}
\Gamma_{(d+1)(d+1)} \subset \Gamma_{[d+3]\cdot 0} \subset \Gamma_{dd}^{(2)}
\end{equation}
\begin{equation} \label{highsyzy}
\Gamma_{[d+3]\cdot 0}/\Gamma_{(d+1)(d+1)} \cong \ker \sym^2 
I^{(d+1)}(d+2) \lra I(2d+4)
\end{equation}
where $I^{(d+1)}(d+2)$ is the space of degree $d+2$ polynomials vanishing
at order $d+1$ along $C_{can}$. Some previous work towards \eqref{serincl}
has been done in \cite{gu3}. Statement \eqref{highsyzy} follows from 
\eqref{serincl}. A generalization of Thm. 1.1 (2) would assert that the 
inclusion
$$ \langle W^{d+1}_{g-1}(C) \rangle \subset \Gamma_{dd}^{(2)} $$
is an isomorphism.

\bigskip
\noindent
{\bf 3.} An important ingredient of the proof of Thm. 1.1 (resp. Thm. 1.2) 
is the use of rank $2$ vector bundles with $3$ (resp. $4$) sections. The 
constructions involved may be viewed as examples of a general construction:
consider, for $n \geq 3$, the subvarieties
$$ \mathcal{W}(n) = \{ [E] \in \modk \: | \: h^0(E) = n \: \text{and} \: E
\: \text{is globally generated} \} .$$
We associate to any $E \in \mathcal{W}(n)$ a commutative diagram
(as in \eqref{grasscon}: replace $V$ by $H^0(E)$)
\begin{equation} \label{grassco2}
\begin{CD}
C @>\gamma>> Gr(2,H^0(E)^*) \\
@VV\varphi_KV    @VVV \\
|K|^* @>\lambda^*>> \mathbb{P}(\wedge^2 H^0(E)^*) = 
\mathbb{P}^{\binom{n}{2} -1}
\end{CD}
\end{equation}
The definitions of the morphisms $\gamma$ and $\lambda^*$ are as in 
\eqref{defgamma} and \eqref{extprod}.

If $n$ is even, $n = 2d+4$ for $d\geq 0$, the Pl\"ucker space
$\Lambda^2 H^0(E)^*$ carries canonically a symmetric multilinear
form
$$ \tilde{\Pf} (\omega_1, \ldots, \omega_{d+2}) = \omega_1 \wedge \ldots \wedge
\omega_{d+2} \in \Lambda^{2d+4} H^0(E)^* \cong \cc $$
which defines a degree $d+2$ polynomial $\Pf \in \sym^{d+2}(\Lambda^2
H^0(E)^*)$ vanishing to order $d+1$ along the Grassmannian 
$Gr(2,H^0(E)^*)$. Notice that $\Pf$ is the Pfaffian if we represent
$\omega \in \Lambda^2 H^0(E)^*$ as an $n \times n$ skew-symmetric matrix.
Therefore we can define for any $E \in \mathcal{W}(2d+4)$ a polynomial
$Q_E = (\lambda^*)^{-1}(\Pf) \in |I^{(d+1)}(d+2)|$. A straighforward
generalization of Prop. 5.12 leads to

\begin{prop}
For all $d \geq 0$, we have a commutative diagram
\begin{equation}
\begin{CD}
\mathcal{W}(2d+4) @>D>> \mathbb{P} \Gamma_{dd}^{(2)} \\
@VVQV @VV\phi_{d+2}^*V \\
|I^{d+1}(d+2)| @>Ver>> \mathbb{P} \sym^2 I^{d+1}(d+2)              
\end{CD}
\end{equation}
\end{prop}
In order to prove surjectivity of $\phi_{d+2}^*$ (for $d>0$) one needs a
better understanding of the vector space $I^{d+1}(d+2)$. Can one construct
naturally a basis of $I^{d+1}(d+2)$ as in the $d=0$ case (Petri's quadrics)?

If $n$ is odd and $n>3$, the geometry one can extract from diagram 
\eqref{grassco2}
seems more complicated, due to the fact that there is no natural equation
attached to the Grassmannian. If $n=3$, the composite $\lambda^* \circ 
\varphi_K :
C \lra \pp^2 = \pp(\Lambda^2 H^0(E)^*)$ is the morphism described in
remark 2.5.

\bigskip
\noindent
{\bf 4.} Finally, we note that the bundles $E_W$, which were introduced 
in section 2.3, can be used to work out a vector bundle-theoretical proof
of Welters' theorem \cite{welt}, i.e. the base locus of $\gg$ equals the 
surface $C-C$ for $g \geq 5$.

\bigskip

\bigskip
\noindent
\flushleft{Christian Pauly \\
Laboratoire J.-A. Dieudonn\'e \\
Universit\'e de Nice Sophia Antipolis \\
Parc Valrose \\
F-06108 Nice Cedex 02 \\ France \\
E-mail: pauly@math.unice.fr}

\bigskip
\noindent
\flushleft{Emma Previato \\
Department of Mathematics \\
Boston University \\
Boston, MA 02215-2411 \\
U.S.A. \\
E-mail: ep@math.bu.edu}

\end{document}